\documentclass[12pt,a4paper]{article}
\usepackage[T2A]{fontenc}
\usepackage[cp1251]{inputenc}
\usepackage[english, russian]{babel}
\usepackage{indentfirst}
\usepackage{amsmath,amssymb,amsthm,amsfonts,amscd}
\usepackage{graphicx}

\newcommand{\Div}{\mathop{\fam0 div}}

\newcommand{\grad}{\mathop{\fam0 grad}}

\newcommand{\Ker}{\mathop{\fam0 Ker}}
\newcommand{\Hom}{\mathop{\fam0 Hom}}

\renewcommand{\le}{\leqslant}
\renewcommand{\ge}{\geqslant}
\def\to{\rightarrow}

\newtheorem{theorem}{Теорема}
\newtheorem{lemma}{Лемма}
\newtheorem{statement}{Утверждение}

\begin{document}
\sloppy

\title{Высокочастотные асимптотики периодических по времени решений систем дифференциальных уравнений в критическом случае}

\author{Левенштам В. Б., Нгуен Л. К., Ишмеев М. Р.}

\maketitle

\noindent{УДК 519.4+513.88}

\begin{abstract}

Для двух линейных эволюционных систем дифференциальных уравнений --- нормальной системы обыкновенных дифференциальных уравнений и системы уравнений в частных производных с оператором Стокса в главной части --- с быстро осцилирующими по времени коэффициентами в случае кратного вырождения усредненного стационарного оператора рассмотрена задача о периодических по времени решениях. Доказаны результаты о существовании и единственности таких решений и построены с обоснованием асимптотические разложения последних. Для нормальных систем обыкновенных дифференциальных уравнений, кроме того, установлена сходимость асимптотического ряда в обычном смысле и исследованы вопросы об устойчивости по Ляпунову решения.

%The time-periodical solutions problem was considered for two evolution differential equations systems with rapidly oscillating by time coefficients in the case of the multiple degeneration of the averaged stationary operator: the normal ordinary differential equations system and the partial derivatives differential equations system with Stokes operator in the main part. The results about the existense and the uniqueness of these solutions were approved and their asymptotic expansions were constructed and substanded. Moreover, the  asymtotic  sums  were approved and Lyapynov stability problems were researched for normal ordinary differential equations system.  

\end{abstract}

\noindent \textbf{Ключевые слова}: линейная нормальная система дифференциальных уравнений, оператор Стокса, высокочастотная осцилляция коэффициентов, метод усреднения, полная асимптотика, критический случай, метод пограничного слоя.

\section*{Введение}

В работах \cite{Do1}, \cite{Do2} рассматривается нормальная система линейных дифференциальных уравнений с быстро осциллирующими коэффициентами, для которой предельная (усредненная в смысле \cite{kop:bib1}) стационарная задача имеет простое нулевое собственное значение. Установлены существование и единственность периодического решения данной системы и построена его полная асимптотика по степеням малой величины, обратной частоте осцилляций коэффициентов системы, в двух случаях: когда собственный вектор, отвечающий простому нулевому собственному значению предельной стационарной задачи, не имеет обобщенных (в смысле Вишика-Люстерника \cite{kop:bib4}) присоединенных векторов \cite{Do1}, и когда имеет обобщенный присоединенный  вектор \cite{Do2}. Исследованы также вопросы устойчивости и неустойчивости по Ляпунову этого решения и доказана сходимость в обычном смысле асимптотического  ряда.

В работах \cite{Gus1}--\cite{Lev4} рассматривались линейные параболические задачи в случае простого вырождения стационарной предельной задачи при отсутствии и наличии обобщенных присоединенных векторов 1-го порядка. Доказаны результаты о существовании и единственности периодических по времени решений, а также при помощи метода пограничного слоя \cite{Vishik2} построены их асимптотические разложения. Аналогичные результаты получены в \cite{Our1}, \cite{Our2} для линейных систем уравнений в частных производных с оператором Стокса в главной части также в случае простого вырождения.

В первом параграфе работы рассматривается возмущённая нормальная система того же вида, что в  \cite{Do1}, \cite{Do2}. Однако предельная стационарная задача имеет нулевое собственное значение произвольной конечной кратности. Предполагается, что все собственные векторы предельной задачи , отвечающие нулевому собственному значению, не имеют обобщенных присоединенных (в смысле Вишика-Люстерника) векторов. Для указанной системы решены те же, что в \cite{Do1}, \cite{Do2}, вопросы.

Во втором параграфе работы исследуется линейная система уравнений в частных производных с оператором Стокса в главной части в случае кратного вырождения стационарной предельной задачи. Доказаны существование и единственность периодического по времени решения, а также при помощи метода погранслоя построена и обоснована его асимптотика в равномерных по области метриках.

Отметим, что результаты \S 1 были получены в магистерской диссертации Л. К. Нгуена в 2012 г. В работе Л. И. Сазонова \cite{Sazon} доказаны существование и единственность периодических решений, а также построены полные асимптотики последних для линейных дифференциальных уравнений с ограниченными быстро осциллириющими коэффициентами в банаховых пространствах в случае произвольного (в частности кратного) вырождения. В работе Л. И. Сазонова \cite{Sazon2} аналогичные \cite{Sazon} результаты получены для уравнений  с неограниченными операторными коэффициентами в банаховых пространствах. Актуальный для приложений вопрос о погранслойных составляющих асимптотик, приближающих решение в равномерных метриках, в столь широкой постановке \cite{Sazon2}, естественно, не расматривался.

Работа частично поддержана Грантом РФФИ ( N12-01-00402-a ).

\section{Нормальная система обыкновенных дифференциальных уравнений}

\subsection{Основной результат}

Обозначим через $M_n(C)$ множество квадратных матриц порядка $n$ с комплексными элементами, а через $C^n$ --- множество $n$-мерных вектор-столбцов с комплексными компонентами. Пусть $m,n\in N$; $A_0,B_l\in M_n(C)$; $d_l\in C^n$, $0\le|l|\le m$. Рассмотрим задачу о $\frac{2\pi}{\omega}$-периодическом решении системы линейных обыкновенных дифференциальных уравнений
\begin{equation}\label{kop:eq1}
\frac{dx}{dt}=\left(A_0+\frac{1}{\omega}B_0\right)x+\sum\limits_{1\le |l|\le m} (B_l x+d_l)e^{il\omega t}+d_0, 
\end{equation}
где $\omega$ --- большой параметр. Пусть $\lambda=0$ --- собственное значение матрицы $A_0$, $\{a_1 ,a_2 ,\dots,a_s\}$ --- некоторый базис собственного подпространства $\ker A_0$, $(1\le s\le n)$. Введём в рассмотрение матрицу
$$
A_1=B_0+ \sum\limits_{1\le |l|\le m} \frac{B_{-l} B_{l}}{il}.
$$
Будем предполагать, что все собственные векторы с нулевым собственным значением матрицы $A_0$ не имеют присоединенных относительно матриц $A_0$, $A_1$ векторов (см., например \cite{kop:bib4}), т.е. уравнение
$$
A_0 x+A_1(\mu_1 a_1+\mu_2 a_2+\dots+\mu_s a_s)=0
$$
не имеет решений для любого нетривиального набора чисел $\mu_1,\dots, \mu_s$. 

При указанных условиях имеет место следующее утверждение

\begin{theorem} Система (\ref{kop:eq1}) при достаточно больших $\omega$ имеет единственное $\frac{2\pi}{\omega}$-периодическое решение, представимое в виде равномерно и абсолютно сходящегося на оси $t\in R$ ряда
\begin{equation}\label{kop:eq2}
x(t)=\omega\sum\limits_{j=1}^s C_{-1}^j a_j+\sum\limits_{k=0}^\infty\frac{1}{\omega^k}\left(x_k+y_k(\omega t)+\sum\limits_{j=1}^s C_k^j a_j\right),
\end{equation}
где $C_{-1}^j, C_k^j\in C$; $x_k\in C^n$; $y_k: R\to C^n$ --- $2\pi$-периодические вектор-функции с нулевым средним, $j=1,\dots,s$, $k=0,1,2,\dots$, которые эффективно определяются. При этом существует такая постоянная $M_0$, что для всех $t\in R$ при достаточно больших $\omega$ имеет место оценка
\begin{equation}\label{kop:eq3}
\left|x(t)-\omega\sum\limits_{j=1}^s C_{-1}^j a_j+\sum\limits_{k=0}^r\frac{1}{\omega^k}\left(x_k+y_k(\omega t)+\sum\limits_{j=1}^s C_k^j a_j\right)\right|\le\frac{M_0^{r+1}}{\omega^{r+1}},
\end{equation}
$r=0,1,2,\dots$, свидетельствующая о том, что ряд (\ref{kop:eq2}) - асимптотический.
\end{theorem}

Напомним при этом, что средним непрерывной $T$-периодической функции $a(\tau)$ называется величина $\langle a(\tau)\rangle=\frac{1}{T}\int\limits_{0}^{T}a(\tau)d\tau$, а среднее вектор-функции $(a_1(\tau),a_2(\tau),\dots,a_n(\tau))$ определяется равенством $\langle(a_1(\tau),a_2(\tau),\dots,a_n(\tau))\rangle =(\langle a_1(\tau)\rangle ,\langle a_2(\tau)\rangle ,\dots,\langle a_n(\tau)\rangle)$.

Под эффективностью нахождения коэффициентов ряда (\ref{kop:eq2}) понимаем тот факт, что нахождение каждого из них сводится к конечному числу арифметических действий над известными числами и функциями (см. п.3). Доказательство теоремы 1 изложено в п.3-5.

Заметим, что решение $x(t)$ будет вещественным, если коэффициенты правой части системы (\ref{kop:eq1}) вещественны.

Предполагается теперь, что матрицы $A_0$, $B_0$ вещественны, а матрицы $B_{-l}$, $B_l$ $(1\le l\le m)$ комплексно сопряжены, т.е. коэффициент при $x$ правой части системы (\ref{kop:eq1}) веществен. Рассмотрим соответствующую (\ref{kop:eq1}) однородную систему
$$
\frac{dx}{dt}=\left(A_0+\frac{1}{\omega}B_0\right)x+\sum\limits_{1\le l\le m} B_l xe^{il\omega t}.
$$
Производя здесь последовательно бесконечное число замен переменных Крылова-Боголюбова (см.\cite{kop:bib1}, с.24) ( $x=y+ \sum\limits_{1\le l\le m} \frac{1}{il\omega} B_l y e^{il\omega t}$ --- первая из них), придем к формальной системе
\begin{equation}\label{kop:eq4}
\frac{dy}{dt}=(A_0+\frac{1}{\omega}A_1+\frac{1}{\omega^2}A_2+\dots)y\equiv A_\omega y.
\end{equation}
Рассмотрим характеристический многочлен $|\lambda E-A_\omega|=\lambda^n+\alpha_1 \lambda^{n-1}+\dots+\alpha_n$ и составим по нему матрицу Гурвица
$$
\Gamma (\omega)=\begin{pmatrix} \alpha_1 & 1 & 0 & 0 & \cdots & 0 \\ \alpha_3 & \alpha_2 & \alpha_1 & 1 & \cdots & 0 \\ \vdots & \vdots & \vdots & \vdots & \ddots & \vdots \\ \alpha_{2n-1} & \alpha_{2n-2} & \alpha_{2n-3} & \alpha_{2n-4} & \cdots & \alpha_n \end{pmatrix},
$$
где $\alpha_s=0$, $\forall s>n$. Детерминантами Гурвица, относящимся к формальной системе дифференциальных уравнений с постоянными коэффициентами (\ref{kop:eq4}), являются все диагональные миноры $D_{j\omega}(j=1,\dots,n)$ матрицы Гурвица $\Gamma(\omega)$. Пусть все эти детерминанты $D_{j\omega}$ $(j=1,\dots,n)$, отличны от нуля. Разложим их в формальные ряды по степеням $\omega^{-1}$
$$
\begin{array}{c}
D_{1\omega}=\alpha_1=\sum\limits_{q=0}^\infty\omega^{-q}d_{1q},\quad D_{2\omega}=\begin{vmatrix}\alpha_1 & 1\\ \alpha_3 & \alpha_2\end{vmatrix}=\sum\limits_{q=0}^\infty\omega^{-q}d_{2q},\dots,\\ D_{n\omega}=\det\Gamma=\sum\limits_{q=0}^\infty\omega^{-q}d_{nq}.
\end{array}
$$
Обозначим через $d_{jq_j}$, $j=1,2,\dots,n$ --- первые ненулевые коэфициенты этих разложений. Имеет место следующее утверждение об устойчивости и неустойчивости по Ляпунову $\frac{2\pi}{\omega}$-периодического решения системы (\ref{kop:eq1}).  
\begin{theorem} Если в разложениях по степеням $\omega^{-1}$ детерминантов Гурвица, относящимся к формальной системе (\ref{kop:eq4}), все первые ненулевые коэффициенты $d_{jq_j}(j=1,\dots,n)$ положительны, то при достаточно больших $\omega$ решение $x(t)$ системы (\ref{kop:eq1}) устойчиво по Ляпунову. Если же хотя бы один из первых ненулевых коэффициентов $d_{jq_j}(j=1,\dots,n)$  отрицателен, то решение $x(t)$ системы (\ref{kop:eq1}) неустойчиво.   
\end{theorem}

Заметим, что если среди детерминантов Гурвица $D_{j\omega}(j=1,\dots,n)$ существуют тождественные нулю (хотя бы один), то об устойчивости и неустойчивости системы (\ref{kop:eq1})  ничего сказать нельзя. 
Действительно, рассмотрим следующие системы:
\begin{equation}\label{kop:eqa}
\frac{dx}{dt}=\left[\begin{pmatrix}0 & 0 & -1 \\ 0 & 0 & 0 \\ 0 & 0 & -1 \end{pmatrix}+\frac{1}{\omega}\begin{pmatrix}1 & 0 & -1 \\ 1 & 1 & -1 \\ 2 & 2 & -2 \end{pmatrix}\right]x\equiv\left(A_0+\frac{1}{\omega}B_{01}\right)x\equiv A_{1\omega}x,
\end{equation}
и
\begin{equation}\label{kop:eqb}
\frac{dx}{dt}=\left[\begin{pmatrix}0 & 0 & -1 \\ 0 & 0 & 0 \\ 0 & 0 & -1 \end{pmatrix}+\frac{1}{\omega}\begin{pmatrix}1 & 0 & -1 \\ 0 & 1 & 0 \\ 2 & 0 & -2 \end{pmatrix}\right]x\equiv\left(A_0+\frac{1}{\omega}B_{02}\right)x\equiv A_{2\omega}x,
\end{equation}
Очевидно, для систем (\ref{kop:eqa}), (\ref{kop:eqb}) выполняются все условия теоремы 1, т.е. матрица $A_0$ имеет двукратное нулевое собственное значение, и все её собственные векторы, отвечающие нулевому собственному значению, не имеют присоединённых относительно $A_0$, $B_{01}$ (и также относительно $A_0$, $B_{02}$ векторов). При этом
$$
|\lambda E - A_{1\omega}|=\lambda^3+\lambda^2+\frac{1}{\omega^2}\lambda+\frac{1}{\omega^2}, \quad |\lambda E - A_{2\omega}|=\lambda^3+\lambda^2-\frac{1}{\omega^2}\lambda-\frac{1}{\omega^2}
$$
Соответствующие детерминанты Гурвица систем (\ref{kop:eqa}), (\ref{kop:eqb})  равны при этом: $D^{(a)}_{1\omega}=D^{(b)}_{1\omega}=1$, $D^{(a)}_{2\omega}=D^{(b)}_{2\omega}=D^{(a)}_{3\omega}=D^{(b)}_{3\omega}=0$. Система (\ref{kop:eqa}), устойчива по Ляпунову,поскольку все корни её характеристического уравнения имеют неположительные вещественные части. Система (\ref{kop:eqb}) неустойчива, так как её характеристическое уравнение имеет корень $\lambda=\omega^{-1}>0$. 

Доказательство теоремы 2 опускается. Отметим лишь, что оно базируется на работе \cite{kop:bib5}, а также использует идеи и подходы классической теории устойчивости (см., например \cite{kop:bib6}).

\subsection{Нахождение коэффициентов асимптотического ряда (\ref{kop:eq2})}

Формальное асимптотическое разложение решения системы (\ref{kop:eq1}) на основании метода двухмасштабных разложений и метода Вишика-Люстерника (см.\cite{kop:bib4}) будем искать в виде (\ref{kop:eq2}).
Подставляя (\ref{kop:eq2}) в (\ref{kop:eq1}), получим
\begin{equation}\label{kop:eq5}
\begin{array}{c}
\omega\frac{dy_0}{d\tau}+\sum\limits_{k=0}^\infty\frac{1}{\omega^k}\frac{dy_{k+1}}{d\tau}=\sum\limits_{k=0}^\infty\frac{A_0}{\omega^k}(x_k+y_k(\tau))+\\
+\frac{B_0}{\omega}\left[\omega\sum\limits_{j=1}^s C_{-1}^ja_j+\sum\limits_{k=0}^\infty\frac{1}{\omega^k}\left(x_k+y_k(\tau)+\sum\limits_{j=1}^s C_{k}^ja_j\right)\right]+\\
+\sum\limits_{1\le |l|\le m}B_l\left[\omega\sum\limits_{j=1}^s C_{-1}^ja_j+\sum\limits_{k=0}^\infty\frac{1}{\omega^k}\left(x_k+y_k(\tau)+\sum\limits_{j=1}^s C_{k}^ja_j\right)\right]e^{ilt}+\\
+\sum\limits_{1\le |l|\le m}d_le^{ilt}+d_0,
\end{array}
\end{equation}
где $\tau=\omega t$. Приравняем коэффициенты в левой и правой частях равенства (\ref{kop:eq5}) при одинаковых степенях $\omega$. Начнём со старшей степени $\omega$:
$$
\frac{dy_0}{d\tau}= \sum\limits_{1\le |l|\le m} B_l\left(\sum\limits_{j=1}^sC_{-1}^ja_j\right)e^{il\tau},\quad \langle y_0(\tau) \rangle = 0 
$$  
так что
\begin{equation}\label{kop:eq6}
y_0(\tau)=\sum\limits_{1\le |l|\le m} \frac{B_l}{il}\left(\sum\limits_{j=1}^s C_{-1}^j a_j\right)e^{il\tau}.
\end{equation}
Равенство коэффициентов в (\ref{kop:eq5}) при $\omega^0$ имеет вид
$$
\begin{array}{c}
\frac{dy_1}{d\tau}=A_0(x_0+y_0) + \sum\limits_{j=1}^sC_{-1}^jB_0a_j + \sum\limits_{1\le |l|\le m} B_l \left(x_0 + y_0 + \sum\limits_{j=1}^sC_0^ja_j\right)e^{il\tau} +\\
+\sum\limits_{1\le |l|\le m} d_le^{il\tau} + d_0. 
\end{array}
$$ 
Подставляя сюда $y_0(\tau)$ из (\ref{kop:eq6}) и применяя к полученному равенству операцию усреднения, получим равенство
$$
0=A_0x_0 + \sum\limits_{j=1}^sC_{-1}^jB_0a_j + \sum\limits_{1\le |l|\le m} \frac{B_{-l}B_l}{il}\left(\sum\limits_{j=1}^sC_{-1}^ja_j\right) + d_0,
$$  
или
\begin{equation}\label{kop:eq7}
A_0x_0 + \sum\limits_{j=1}^sC_{-1}^jA_1a_j + d_0=0 
\end{equation}
Заметим, что число $\lambda=0$ является собственным значением сопряжённой к $A_0$ матрицы $A_0^*$, ему также отвечает $s$ линейно независимых собственных векторов $z_1,z_2,\dots,z_s$ матрицы $A_0^*$. При этом линейная оболочка, натянутая на векторы $z_1,z_2,\dots,z_s$, совпадает со множеством всех собственных векторов матрицы $A_0^*$, отвечающих собственному значению $\lambda=0$. Уравнение (\ref{kop:eq7}) разрешимо относительно $x_0$ тогда и только тогда, когда вектор $\sum\limits_{j=1}^sC_{-1}^jA_1a_j+d_0$ ортогонален этой оболочке, т.е. $\left(\sum\limits_{j=1}^sC_{-1}^jA_1a_j+d_0,z_k\right)=0$, $k=1,\dots,s$.
Последнее равносильно СЛАУ
\begin{equation}\label{kop:eq8}
\begin{pmatrix}(A_1a_1,z_1) & (A_1a_2,z_1) & \cdots & (A_1a_s,z_1) \\ (A_1a_1,z_2) & (A_1a_2,z_2) & \cdots & (A_1a_s,z_2) \\ \vdots & \vdots & \ddots & \vdots \\ (A_1a_1,z_s) & (A_1a_2,z_s) & \cdots & (A_1a_s,z_s) \end{pmatrix}\begin{pmatrix}C^1_{-1} \\ C^2_{-1} \\ \cdots \\ C^s_{-1}\end{pmatrix}=-\begin{pmatrix}(d_0,z_1) \\ (d_0,z_2) \\ \cdots \\ (d_0,z_s)\end{pmatrix}\equiv b_{-1}
\end{equation}
Докажем, что матрица $\Delta=(\sigma_{kj})_{k,j=1}^s\equiv \left( (A_1a_j,z_k)\right)_{k,j=1}^s$ невырождена. Действительно, в противном случае обращалась бы в нуль некоторая линейная комбинация столбцов
$$
\sum\limits_{j=1}^s\mu_j(A_1a_j,z_k)=0,\quad k=1,\dots,s, \quad \sum\limits_{j=1}^s\mu_j^2 \neq 0.
$$   
Но отсюда у собственного вектора $\sum\limits_{j=1}^s\mu_ja_j$ матрицы $A_0$ будет какой-то присоединённый относительно матриц $A_0, A_1$ вектор, что невозможно в силу нашего предложения. Итак матрица $\Delta$ невырождена. Из (\ref{kop:eq8}) теперь находим 
\begin{equation}\label{kop:eq9}
C_{-1}^j=\left(\Delta^{-1}\right)_jb_{-1},\quad j=1,\dots,s
\end{equation}
где $\left(\Delta^{-1}\right)_j-j$ -ая строка обратной к $\Delta$ матрицы $\Delta^{-1}$. Из (\ref{kop:eq7}) находим ещё 
\begin{equation}\label{kop:eq10}
x_0=-W\left(\sum\limits_{j=1}^sC_{-1}^jA_1a_j+d_0\right),
\end{equation}
где $W$ --- обратный оператор (матрица) к сужению оператора (матрицы) $A_0$ на ортогональное дополнение векторов $a_1,a_2,\dots,a_s$. Подставляя $C_{-1}^j$ из (\ref{kop:eq9}) в (\ref{kop:eq6}), определим $y_0(\tau)$. Далее имеем
$$
\begin{array}{c}
\frac{dy_1}{d\tau}=\sum\limits_{1\le |l|\le m} B_l \left(x_0 +  \sum\limits_{j=1}^sC_0^ja_j\right)e^{il\tau} + \sum\limits_{1\le |l|\le m} d_le^{il\tau}+\\
+\sum\limits_{1\le |l|\le m}\frac{A_0B_l}{il}\left(\sum\limits_{j=1}^sC_{-1}^ja_j\right)e^{il\tau}+\\
+\sum\limits_{1\le|l_1|,|l_2|\le m,l_1+l_2\neq 0}\frac{B_{l_2}B_{l_1}}{il_1}\left(\sum\limits_{j=1}^sC_{-1}^ja_j\right)e^{i(l_1+l_2)\tau}, \left(y_0(\tau)\right)=0, 
\end{array}
$$
так что
\begin{equation}\label{kop:eq11}
\begin{array}{c}
y_1(\tau)=\sum\limits_{1\le |l|\le m}\frac{B_l}{il}\left(\sum\limits_{j=1}^sC_0^ja_j\right)e^{il\tau}+\sum\limits_{1\le |l|\le m}\frac{B_lx_0}{il}e^{il\tau}-\\
-\sum\limits_{1\le|l_1|,|l_2|\le m,l_1+l_2\ne 0}\frac{B_{l_2}B_{l_1}}{(l_1+l_2)l_1}\left(\sum\limits_{j=1}^sC_{-1}^ja_j\right)e^{i(l_1+l_2)\tau}-\\
-\sum\limits_{1\le |l|\le m}\frac{A_0B_l}{l^2}\left(\sum\limits_{j=1}^sC_{-1}^ja_j\right)e^{il\tau}+\sum\limits_{1\le |l|\le m}\frac{d_l}{il}e^{il\tau}\equiv\\
\equiv\sum\limits_{1\le |l|\le m}\frac{B_l}{il}\left(\sum\limits_{j=1}^sC_0^ja_j\right)e^{il\tau}+\sum\limits_{1\le |l|\le 2m}\beta_l^1e^{il\tau},
\end{array}
\end{equation}
где $\beta_l^1$, $1\le |l| \le 2m$ --- известные векторы.

Итак, определены числа $C_{-1}^j$, $j=1,\dots,s$, вектор $x_0$ и вектор-функция $y_0(\tau)$. Предположим теперь, что для некоторого целого $p\ge 1$ определены числа $C_{-1}^j, C_0^j,\dots,C_{p-2}^j$ $(j=1,\dots,s)$, векторы $x_0,\dots,x_{p-1}$, вектор-функции $y_0,\dots,y_{p-1}$. Кроме того, вектор-функции $y_r(\tau)(0\le r \le p)$ представимы в виде:
\begin{equation}\label{kop:eq12}
y_r(\tau)= \sum\limits_{1\le |l|\le m}\frac{B_l}{il}\left(\sum\limits_{j=1}^sC_{r-1}^ja_j\right)e^{il\tau}+\sum\limits_{1\le |l|\le (r+1)m}\beta_l^re^{il\tau},
\end{equation}
где $\beta_l^0=0,\beta_l^r(1\le r \le p, 1\le |l| \le (r+1)m$ --- известные векторы. Докажем, что отсюда можно определить числа $C_{p-1}^j(j=1,\dots,s)$, вектор $x_p$ и вектор-функцию $y_p(\tau)$. Докажем также, что
$$
\begin{array}{c}
y_{p+1}(\tau)=\sum\limits_{1\le |l|\le m}\frac{B_l}{il}\left(\sum\limits_{j=1}^sC_p^ja_j\right)e^{il\tau}+\\
+\sum\limits_{1\le |l|\le (p+2)m}\beta_l^{p+1}e^{il\tau},
\end{array}
$$
где $\beta_l^{p+1}, 1\le |l|\le (p+2)m$ --- известные векторы. Действительно, приравняем коэффициенты в левой и правой частях равенства (\ref{kop:eq5}) при $\omega^{-p}$:
$$
\begin{array}{c}
\frac{dy_{p+1}}{d\tau}=A_0(x_p+y_p)+B_0\left(x_{p-1}+y_{p-1}+\sum\limits_{j=1}^sC_{p-1}^ja_j\right)+\\
+\sum\limits_{1\le |l|\le m}B_l\left(x_p+y_p+\sum\limits_{j=1}^sC_p^ja_j\right)e^{il\tau}.
\end{array}
$$ 
Подставляя сюда $y_{p-1}(\tau), y_p(\tau)$ вида (\ref{kop:eq12}), и применяя к полученному равенству операцию усреднения, придем к равенству
\begin{equation}\label{kop:eq13}
A_0x_p+\sum\limits_{j=1}^sC_{p-1}^jA_1a_j+B_0x_{p-1}+\sum\limits_{1\le |l|\le m}B_{-l}\beta_l^p=0.
\end{equation} 
Последнее разрешимо относительно $x_p$ тогда и только тогда, когда выполняются равенства
$$
\left(\sum\limits_{j=1}^sC_{p-1}^jA_1a_j+B_0x_{p-1}+\sum\limits_{1\le |l|\le m}B_{-l}\beta_l^p,z_k\right)=0,\quad k=1,\dots,s,
$$
т.е. 
$$
\Delta\large\left(\begin{array}{c}C_{p-1}^1\\ C_{p-1}^2\\ \cdots \\ C_{p-1}^s\end{array}\right)=b_{p-1},
$$
где $b_{p-1}=(b_{p-1}^1,\cdots,b_{p-1}^s)$, $b_{p-1}^k=\left(B_0x_{p-1}+\sum\limits_{1\le |l|\le m}B_{-l}\beta_l^p,z_k\right)$, $k=1,\dots,s$. Отсюда находим
\begin{equation}\label{kop:eq14}
C_{p-1}^j=(\Delta^{-1})_jb_{p-1}, j=1,\dots,s.
\end{equation}
Подставляя $C_{p-1}^j, j=1,\dots,s$ в (\ref{kop:eq12}), определим вектор-функцию $y_p(\tau)$. Из (\ref{kop:eq13}) находим ещё 
\begin{equation}\label{kop:eq15}
x_p=-W\left(\sum\limits_{j=1}^sC_{p-1}^jA_1a_j+B_0x_{p-1}+\sum\limits_{1\le |l|\le m}B_{-l}\beta_l^p\right).
\end{equation}
Далее
$$
\begin{array}{c}
\frac{dy_{p+1}}{d\tau}=\sum\limits_{1\le |l|\le m}B_l\left(x_p+\sum\limits_{j=1}^sC_p^ja_j\right)e^{il\tau}+\sum\limits_{1\le |l|\le m}\frac{A_0B_l}{il}\left(\sum\limits_{j=1}^sC_{p-1}^ja_j\right)e^{il\tau}+\\
+\sum\limits_{1\le |l|\le (p+1)m}A_0\beta_l^pe^{il\tau}+\\
+\sum\limits_{1\le |l|\le m}\frac{B_0B_l}{il}\left(\sum\limits_{j=1}^sC_{p-2}^ja_j\right)e^{il\tau}+\sum\limits_{1\le |l|\le pm}B_0\beta_l^{p-1}e^{il\tau}+\\
+\sum\limits_{1\le |l_1|,|l_2|\le m, l_1+l_2\neq 0}\frac{B_{l_2}B_{l_1}}{il_1}\left(\sum\limits_{j=1}^sC_{p-1}^ja_j\right)e^{i(l_1+l_2)\tau}+\\
+\sum\limits_{1\le |l_1|\le (p+1)m,1\le|l_2|\le m,l_1+l_2\neq 0}B_{l_2}\beta_{l_1}^pe^{i(l_1+l_2)\tau},\quad \left\langle y_0(\tau)\right\rangle=0,
\end{array}
$$
так что 
\begin{equation}\label{kop:eq16}
y_{p+1}(\tau)=\sum\limits_{1\le |l|\le m}\frac{B_l}{il}\left(\sum\limits_{j=1}^sC_p^ja_j\right)e^{il\tau}+\sum\limits_{1\le |l|\le (p+2)m}\beta_l^{p+1}e^{il\tau},
\end{equation}
где $\beta_l^{p+1}$, $1\le|l|\le (p+2)m$ --- известные векторы, определённые равенством
\begin{equation}\label{kop:eq17}
\begin{array}{c}
\sum\limits_{1\le |l|\le (p+2)m}\beta_l^{p+1}e^{il\tau}=\sum\limits_{1\le |l|\le m}\frac{B_lx_p}{il}e^{il\tau}-\sum\limits_{1\le |l|\le m}\frac{A_0B_l}{l^2}\left(\sum\limits_{j=1}^sC_{p-1}^ja_j\right)e^{il\tau}+\\
+\sum\limits_{1\le |l|\le (p+1)m}\frac{A_0\beta_l^p}{il}e^{il\tau}-\sum\limits_{1\le |l|\le m}\frac{B_0B_l}{l^2}\left(\sum\limits_{j=1}^sC_{p-2}^ja_j\right)e^{il\tau}+\sum\limits_{1\le |l|\le pm}\frac{B_0\beta_l^{p-1}}{l^2}e^{il\tau}+\\
+\sum\limits_{1\le |l_1|,|l_2|\le m,l_1+l_2\neq 0}\frac{B_{l_2} B_{l_1}}{(l_1+l_2)l_1}\left(\sum\limits_{j=1}^sC_{p-1}^ja_j\right)e^{i(l_1+l_2)\tau}+\\
+\sum\limits_{1\le |l_1|\le (p+1)m 1\le |l_2|\le m,l_1+l_2\ne 0}\frac{B_{l_2}\beta_{l_1}^p}{i(l_1+l_2)}e^{i(l_1+l_2)\tau}.
\end{array}
\end{equation}
Таким образом, можно найти все коэффициенты, векторы и вектор-функции, фигурирующие в формальной асимптотике (\ref{kop:eq2}).

\subsection{Сходимость ряда (\ref{kop:eq2}) к решению системы (\ref{kop:eq1})}

Докажем, что существует такое число $\omega_0>0$, что при $\omega>\omega_0$ ряд (\ref{kop:eq2}) с коэффициентами, определёнными в пункте 3, равномерно сходится и действительно является решением системы (\ref{kop:eq1}). В $C^n$ для любых двух векторов $x=(x_1,\dots,x_n)$ и $y=(y_1,\dots,y_n)$ определим скалярное произведение $(x,y)=\sum\limits_{i=1}^nx_i\overline{y_i}$, а также введем норму вектора $| x|=\sqrt{(x,x)}$. Введём в $M_n(C)$ норму матрицы $A=(a_{ij})$:
$$
\parallel A\parallel =\left(\sum\limits_{1\le i, j\le n}| a_{ij}| ^2\right)^{\frac{1}{2}}.
$$ 
Напомним при этом, что для любого вектора $x\in C^n$ и любой матрицы $A\in M_n(C)$ выполняется неравенство $| Ax| \le \parallel A\parallel | x|$.

Не теряя общности, будем считать, что при всех $l:0\le | l| \le m$ выполняются соотношения $\parallel A_0 \parallel \le 1, \parallel B_l \parallel \le 1, | d_l| \le 1$. Действительно, этого всегда можно добиться заменой независимой переменной $t\to pt$, $t\in R$.Также можно считать, что $| a_i|=1, i=1,\dots,s$. Введём ряд обозначений
$$
\theta_0=d_0, \theta_p=B_0x_{p-1}+\sum\limits_{1\le |l| \le m}B_{-l}\beta_{l}^p,
$$
$$
\mu_0=0,\quad \mu_p=\sum\limits_{1\le|l| \le (p+1)m}| \beta_l^p|,\quad p=1,2,\cdots,
$$
где векторы $\beta_l^1(1\le|l| \le 2m)$ определены соотношением (\ref{kop:eq11}), а векторы $\beta_l^p$ ($p\ge 2$, $1\le |l| \le (p+1)m$) --- соотношением (\ref{kop:eq17}). Имеем для коэффициентов ряда (\ref{kop:eq2}) следующие соотношения
$$
\begin{array}{c}
x_p=-W\left(\sum\limits_{j=1}^sC_{p-1}^jA_1a_j+\theta_p\right),\\
y_p(\tau)=\sum\limits_{1\le |l| \le m}\frac{B_l}{il}\left(\sum\limits_{j=1}^sC_{p-1}^ja_j\right)e^{il\tau} + \sum\limits_{1\le |l| \le (p+1)m}\beta_l^pe^{il\tau},
\end{array}
$$
\begin{equation}\label{kop:eq18}
\large\left(\begin{array}{c}C_{p-1}^1\\C_{p-1}^2 \\ \cdots \\ C_{p-1}^s\end{array}\right)=\Delta^{-1}\left(\begin{array}{c}-(\theta_p,z_1) \\ -(\theta_p,z_2) \\ \cdots \\ -(\theta_p,z_s) \\ \end{array}\right),\quad p=0,1,2,\dots
\end{equation}
Из соотношений (\ref{kop:eq18}) очевидно, что существует положительное $L$ такое, что выполняются равенства
\begin{equation}\label{kop:eq19} 
| C_{p-1}^j| \le \frac{L}{s}| \theta_p |,\quad j=1,\dots, s, \quad | x_p| \le L| \theta_p | \forall p=0,1,2,\dots
\end{equation}
В силу (\ref{kop:eq11}) имеем теперь
\begin{equation}\label{kop:eq20}
\begin{array}{c} 
\mu_1=\sum\limits_{1\le |l| \le 2m}| \beta_l^1 | \le \sum\limits_{1\le |l| \le m}\left|\frac{B_lx_0}{il}\right|+\sum\limits_{1\le|l_1|, |l_2|\le m, l_1+l_2\ne 0}\left|\frac{B_{l_2}B_{l_1}}{(l_1+l_2)l_1}\left(\sum\limits_{j=1}^sC_{-1}^ja_j\right)\right| +\\
+ \sum\limits_{1\le |l| \le m}\left|\frac{A_0B_l}{l^2}\left(\sum\limits_{j=1}^sC_{-1}^ja_j\right)\right| + \sum\limits_{1\le |l| \le m}\left|\frac{d_l}{il}\right| \le 2m|x_0| +\\
+(2m)^2L| \theta_0| +2mL| \theta_0 | + 2m\le \left[(2m)^2+4m\right]L| \theta_0| +2m.
\end{array}
\end{equation}
Далее, в силу (\ref{kop:eq17}) для всех $p\ge 2$ выполняется неравенство 
\begin{equation}\label{kop:eq21}
\begin{array}{c}
\mu_p = \sum\limits_{1\le|l|\le (p+1)m}|\beta_l^p|\le\sum\limits_{1\le |l|\le m}\left|\frac{B_lx_{p-1}}{il}\right|+\sum\limits_{1\le |l|\le m}\left|\frac{A_0B_l}{l^2}\left(\sum\limits_{j=1}^sC_{p-2}^ja_j\right)\right|+\\
+\sum\limits_{1\le |l|\le pm}\left|\frac{A_0B_l^{p-1}}{il}\right|+\sum\limits_{1\le |l|\le m}\left|\frac{B_0B_l}{l^2}\left(\sum\limits_{j=1}^sC_{p-3}^ja_j\right)\right|+\sum\limits_{1\le |l|\le (p-1)m}\left|\frac{B_0B_l^{p-2}}{il}\right|+\\
+\sum\limits_{1\le |l_1|,|l_2|\le m, l_1+l_2\neq 0}\left|\frac{B_{l_1}B_{l_2}}{(l_1+l_2)l_1}\left(\sum\limits_{j=1}^sC_{p-2}^ja_j\right)\right|+\\
+\sum\limits_{1\le |l_1|\le pm, 1\le |l_2|\le m, l_1+l_2\neq 0}\left|\frac{B_{l_2}\beta_{l_1^{p-2}}}{il}\right|\le\\
\le 2m|x_{p-1}|+\mu_{p-1}+\mu_{p-2} + 2m\mu_{p-1}+2mL|\theta_{p-1}|+2mL|\theta_{p-2}|+\\
+(2m)^2L|\theta_{p-1}|\le(2m+1)\mu_{p-1}+\mu_{p-2}+[(2m)^2+4m]L|\theta_{p-1}||+2mL|\theta_{p-2}|.
\end{array}
\end{equation}
Положим $K=2m+2$, $\varphi_0=2m$, $\varphi_p=(2m+1)\mu_p+\mu_{p-1}+2mL| \theta_{p-1}|, p=1,2,\dots$ Докажем следующий результат для величин $\varphi_p, p\ge 1$.
\begin{lemma}
При всех $p\ge 1$ имеет место неравенство
$$
\varphi_p\le K^3L| \theta_{p-1} |+K\varphi_{p-1},
$$
\end{lemma}
Действительно, в силу (\ref{kop:eq20}) имеем
$$
\begin{array}{c}
\varphi_1=(2m+1)\mu_1+2mL| \theta_0 | \le (2m+1)\left\{\left[(2m)^2+4m\right]L| \theta_0 | +2m \right\} +\\
+2mL| \theta_0 | \le (2m+2)^3L| \theta_0 |+(2m+2)2m=K^3L| \theta_0 | +k\varphi_0.
\end{array}
$$ 
При $p\ge 2$ в силу (\ref{kop:eq20}), (\ref{kop:eq21}) имеют место неравенства
$$
\begin{array}{c} 
\varphi_p=(2m+1)\mu_p+\mu_{p-1}+2mL| \theta_{p-1}| \le\\
\le(2m+1)\left((2m+1)\mu_{p-1}+\mu_{p-2}+\left[(2m)^2+4m\right]L| \theta_{p-1} |\right. +\\
\left.+2mL| \theta_{p-2} | \right)+\mu_{p-1}+2mL| \theta_{p-1} |=\\ =\left[(2m+1)^2+1\right]\mu_{p-1}+(2m+1)\mu_{p-2}+\\
+2m(2m+1)L| \theta_{p-2} | +\left((2m+1)\left[(2m)^2+4m\right]+2m\right)L| \theta_{p-1} \le\\
\le (2m+2)^3L| \theta_{p-1} | +(2m+2)\varphi_{p-1}=K^3L| \theta_{p-1}| +K\varphi_{p-1}.
\end{array}
$$  
Лемма 1 доказана.

Докажем теперь следующий результат для векторов $\theta_p$, $p=0,1,2, \dots$
\begin{lemma}
Для любого целого числа $p\ge 0$ выполняется неравенство
$| \theta_p | \le K^p(KL+1)^p$.
\end{lemma}
Доказательство. Имеем $| \theta_0 | = | d_0 | \le 1$. В силу соотношения (\ref{kop:eq20}) находим
$$
\begin{array}{c} 
| \theta_1 | = \left| B_0x_0 + \sum\limits_{1\le |l| \le m}B_{-l}\beta_l^1 \right| \le | x_0 | + \mu_1 \le \\
\le \left[(2m)^2+4m+1\right]L| \theta_0 | +2m < K^2L|  \theta_0 | +2m < K^2L+K.
\end{array}
$$ 
Используя соотношение (\ref{kop:eq21}) и лемму 1, для любого $p \le 2$ имеем
$$
\begin{array}{c}   
| \theta_p | = \left| B_0x_{p-1} + \sum\limits_{1\le |l| \le m}B_{-l}\beta_l^p \right| \le | x_{p-1}| +\mu_p \le\\
\le (2m+1)\mu_{p-1}+\mu_{p-2}+\left[(2m)^2+4m+1\right]L| \theta_{p-1} | +2mL| \theta_{p-2} | \le\\
\le (2m+1)^2L| \theta_{p-1} | +\varphi_{p-1}\le K^2L| \theta_{p-1} | +K^3L| \theta_{p-2} | +K\varphi_{p-2}\le\\
\le K^2L| \theta_{p-1} | +K^3L| \theta_{p-2} | +K^4L| \theta_{p-3} | +K^2\varphi_{p-3}\le \dots \le\\
\le K^2L| \theta_{p-1}| +K^3L| \theta_{p-2}| +\dots +K^pL| \theta_1 | +K^{p-2}\varphi_1 \le\\
\le K^2L| \theta_{p-1}| +K^3L| \theta_{p-2}| +\dots +K^{p+1}L| \theta_0 | +2mK^{p-1}=\\
=\sum\limits_{j=1}^pK^{j+1}L| \theta_{p-j}| +2mK^{p-1}\le\\
\le K^2L\left(\sum\limits_{j=1}^{p-1}K^{j+1}L| \theta_{p-1-j}| +2mK^{p-2}\right)+\sum\limits_{j=2}^pK^{j+1}L| \theta_{p-j}| +2mK^{p-1}=\\
=K(KL+1)\left(\sum\limits_{j=1}^{p-1}K^{j+1}L| \theta_{p-1-j}| +2mK^{p-2}\right)\le\\
\le K^2(KL+1)^2\left(\sum\limits_{j=1}^{p-2}K^{j+1}L| \theta_{p-2-j}| +2mK^{p-3}\right)\le \dots \le\\
\le K^{p-1}(KL+1)^{p-1}(K^2L| \theta_0 | +2m) < K^p(KL+1)^p. 
\end{array}
$$ 
Лемма 2 доказана. В силу доказательства леммы 2 справедливо также неравенство
\begin{equation}\label{kop:eq22}
\mu_p \le K^p(KL+1)^p, \quad p \ge 1
\end{equation}

Оценим теперь норму общего члена ряда (\ref{kop:eq2}). В силу леммы 1 и соотношения (\ref{kop:eq22}) для любого $p\ge 0$, имеем 
$$
\begin{array}{c}   
\left| x_p+y_p(\omega t)+\sum\limits_{j=1}^sC_p^ja_j\right| \le |x_p|+\left| \sum\limits_{1\le |l| \le m}\frac{B_l}{il}\left(\sum\limits_{j=1}^sC_{p-1}^ja_j\right)e^{il\tau}+\right.\\
\left.+\sum\limits_{1\le |l|\le(p+1)m}\beta_l^pe^{il\tau}\right| + \left| \sum\limits_{j=1}^sC_p^ja_j\right| \le\\
\le L| \theta_p | +2mL| \theta_p | +\mu_p + L| \theta_{p+1}| = (2mL+1)| \theta_p | +\mu_p +L| \theta_{p+1}| \le\\
\le (2mL+L+1)K^p(KL+1)^p+LK^{p+1}(KL+1)^{p+1}=CK^p(KL+1)^p,
\end{array}
$$  
где $C=(2mL+L+1)+LK(KL+1)$ --- некоторое положительное постоянное. Положим $\omega_0 = K(KL+1)$. Таким образом, при достаточно больших $\omega > \omega_0$ ряд (\ref{kop:eq2}) сходится абсолютно и равномерно. Отсюда и в силу последнего соотношения нетрудно вывести оценку (\ref{kop:eq3}).

Проверим теперь, что сумма ряда (\ref{kop:eq2}) при достаточно больших $\omega$ действительно является решением системы (\ref{kop:eq1}). Пусть
$$
S_{\omega}^{(r)}=\omega \sum\limits_{j=1}^sC_{-1}^ja_j+\sum\limits_{k=0}^r \frac{1}{\omega^k}\left(x_k + y_k(\omega t)+\sum\limits_{j=1}^sC_k^ja_j\right)
$$   
$r$-ая частичная сумма ряда (\ref{kop:eq3}). Тогда $S_{\omega}^{(r)}$ при $\omega > \omega_0$ равномерно относительно $t \in R$ сходится к некоторой вектор-функции $x(t) : S_{\omega}^{(r)}\rightrightarrows x(t)$. Из равенства 
$$
\begin{array}{c} 
\frac{dS_{\omega}^{(r+2)}(t)}{dt}=\left(A_0 + \sum\limits_{1\le |l| \le m}B_l e^{il\omega t}\right)S_{\omega}^{(r+1)}(t)
+\frac{1}{\omega}B_0S_{\omega}^{(r)}(t)+\sum\limits_{0\le |l| \le m}d_le^{il\omega t},\\
r=1,2,\dots,
\end{array}
$$
следует равномерная сходимость
$$
\begin{array}{c} 
\frac{dS_{\omega}^{(r)}(t)}{dt}=\left(A_0 + \frac{1}{\omega}B_0\right)x(t)+\sum\limits_{1\le |l| \le m}B_l x(t)e^{il\omega t}
+\sum\limits_{0\le |l|\le m}d_le^{il\omega t} при r \to \infty.
\end{array}
$$
Отсюда находим
$$
\begin{array}{c} 
S_{\omega}^{(r)}(t)-S_{\omega}^{(r)}(0)\rightrightarrows\\
\rightrightarrows \int\limits_{0}^{t}\left[\left(A_0 + \frac{1}{\omega}B_0\right)x(\xi)+
\sum\limits_{1\le |l| \le m}B_l x(\xi)e^{il\omega \xi}+\sum\limits_{0\le |l| \le m}d_le^{il\omega \xi}\right]d\xi, 
\end{array}
$$
так что
$$
\begin{array}{c} 
x(t)-x(0)=\lim\limits_{r \to \infty}\left(S_{\omega}^{(r)}(t)-S_{\omega}^{(r)}(0)\right)=\\
=\int\limits_{0}^{t}\left[\left(A_0 + \frac{1}{\omega}B_0\right)x(\xi)+\sum\limits_{1\le |l| \le m}B_l x(\xi)e^{il\omega \xi}+\sum\limits_{0\le |l| \le m}d_le^{il\omega \xi}\right]d\xi.
\end{array}
$$
Дифференцируя обе части последнего равенства, получим
$$
\frac{dx(t)}{dt}=\left(A_0 + \frac{1}{\omega}B_0\right)x(t)+\sum\limits_{1\le |l| \le m}(B_lx(t)+d_l)e^{il\omega t}+d_0,
$$
т.е. $x(t)$ --- решение системы (\ref{kop:eq1}). 

\subsection{Единственность $\frac{2\pi}{\omega}$-периодического по $t$ решения системы (\ref{kop:eq1})}

Докажем, что система (\ref{kop:eq1}) при указанных условиях и достаточно больших $\omega$ имеет не более одного $\frac{2\pi}{\omega}$-периодического решения. Ввиду $\frac{2\pi}{\omega}$-периодичности по $t$ коэффициентов системы (\ref{kop:eq1}) с каждым её решением $x(t)$ решением является $y(t)=x\left(t+\frac{2\pi}{\omega}\right)$. Поэтому нам достаточно доказать, что система (\ref{kop:eq1}) имеет не более одного $N\frac{2\pi}{\omega}$-периодического решения для некоторого натурального $N$ при достаточно больших $\omega$.

Пусть $\varepsilon >0$, $ A_1=(b_{ij})$, $\quad A_{\varepsilon}=A_0 + \varepsilon A_1$. Не нарушая общности, будем считать, что матрица $A_0$ задана в жордановой форме:
$$
A_0=\begin{pmatrix}
K_1 & \cdots & 0 & 0 & 0 & \cdots & 0 \\
\vdots & \ddots & \vdots & \vdots & \vdots & \ddots & \vdots \\
0 & \cdots & K_s & 0 & 0 & \cdots & 0 \\
0 & \cdots & 0 & \lambda_{p+1} & \beta_{p+1} & \cdots & 0 \\
0 & \cdots & 0 & 0 & \lambda_{p+2} & \cdots & 0 \\
\vdots & \ddots & \vdots & \vdots & \vdots & \ddots & \vdots \\
 0 & \cdots & 0 & 0 & 0 & \cdots & \beta_{n-1}\\
 0 & \cdots & 0 & 0 & 0 & \cdots & \lambda_{n}
\end{pmatrix},
$$
где 
$$
K_r=\begin{pmatrix}
0 & 1 & \cdots & 0 \\
0 & 0 & \cdots & 0 \\
\vdots & \vdots & \ddots & \vdots \\
0 & 0 & \cdots & 1 \\
0 & 0 & \cdots & 0
\end{pmatrix}
$$
--- клетка Жордана порядка $e_r$ $(r=1,\dots ,s)$, $\sum\limits_{r=1}^s e_r =p$ --- кратность собственного значения $\lambda = 0$ матрицы $A_0$, $\lambda_{p+1}\lambda_{p+2}\cdots \lambda_n \neq 0$, а $\beta_j(j=p+1,\dots ,n-1)$ равно 1 или 0.

Положим $p_r = \sum\limits_{j=1}^r e_j+1$ $(r=1,\dots, s)$. Тогда
$$
a_1=(10\dots 0)^{T}, \quad a_2=(0\dots 010\dots 0)^{T},\dots,a_s=(0\dots 010\dots 0)^{T},
$$
где для каждого вектора $a_r \quad (r=2,\dots , s)$ элемент 1 находится в $p_{r-1}$-ой строке соответственно, являются линейно независимыми собственными векторами матрицы $A_0$, отвечающими нулевому собственному значению, причём любой собственный вектор $a$, отвечающий нулевым собственному значению матрицы $A_0$ имеет вид
\begin{equation}\label{kop:eq23}
a=\mu_1a_1+\mu_2a_2+\dots +\mu_sa_s=(\mu_10\dots 0\mu_20\dots 0\mu_s0\dots 0)^{T}, \quad \sum\limits_{j=1}^s|\mu_j|\neq 0
\end{equation}
где компоненты $\mu_r \quad (r=1,\dots ,s)$ вектора $a$ находятся в $p_{r-1}$-ой строке соответственно. Так как все собственные векторы матрицы $A_0$, отвечающие собственному значению $\lambda=0$, не имеют присоединённых относительно матриц $A_0$, $A_1$ векторов, то уравнение $A_0x=-A_1a$ не разрешимо для всех векторов $a$ вида (\ref{kop:eq23}). Поэтому для любого набора чисел $\mu_1,\mu_2,\dots,\mu_s$, таких что $\sum\limits_{j=1}^s|\mu_j|\ne 0$, выполняется хотя бы одно из следующих соотношений
$$  
\begin{array}{c} 
b_{q_11}\mu_1+b_{q_1p_1}\mu_2+\dots+b_{q_1p_{s-1}}\mu_s\ne 0,\\
\\ b_{q_21}\mu_1+b_{q_2p_1}\mu_2+\dots+b_{q_2p_{s-1}}\mu_s\ne 0,\dots,\\
\\ b_{q_s1}\mu_1+b_{q_sp_1}\mu_2+\dots+b_{q_sp_{s-1}}\mu_s\ne 0, \quad q_r=\sum\limits_{j=1}^re_j.
\end{array}
$$
Следовательно, определитель
\begin{equation}\label{kop:eq24}
\begin{vmatrix}
b_{q_11} & b_{q_1p_1} & \cdots & b_{q_1p_{s-1}}\\b_{q_21} & b_{q_2p_1} & \cdots & b_{q_2p_{s-1}}\\ \vdots & \vdots & \ddots & \vdots \\ b_{q_s1} & b_{q_sp_1} & \dots & b_{q_sp_{s-1}}
\end{vmatrix}
\neq 0.
\end{equation}

Согласно регулярной теории возмущений (см., например \cite{kop:bib7}) существует нумерация собственных значений $\lambda_j$, $\lambda_j^{\varepsilon}$, $j=1,\dots, n$ матриц $A_0$, $A_{\varepsilon}$ соответственно такая что $\lambda_j^{\varepsilon}=\lambda_j+o (1)$ при $\varepsilon \to 0$, $j=1,\dots, n$, причём $\lambda_j=0$, $j=1,\dots, p$. Кроме того, как известно, все собственные значения $\lambda_j^{\varepsilon}, \quad j=1,\dots, p$, матрицы $A_{\varepsilon}$ при достаточно малых $\varepsilon$ разлагаются в ряд по целым или дробным степеням $\varepsilon$. Пусть $p_j>0$ (целые или дробные, $ j=1,\dots, p$) --- наименьшие степени $\varepsilon$ этих разложений с соответствующими первыми ненулевыми коэффициентами $\varphi_j$ $(j=1,\dots, p)$, т. е.
$$
\lambda_j^{\varepsilon}=\varphi_j\varepsilon^{p_j}+o(\varepsilon^{p_j}),\quad \varepsilon \to 0, \quad \varphi_j\ne 0, \quad j=1,\dots, p,
$$
Докажем от противного, что $0<p_j\le1$ для всех $1\le j\le p$. Действительно, в противном случае положим, что $\lambda_0^{\varepsilon}=\varphi_0\varepsilon^{p_0}+\sum\limits_{i=1}^{\infty}\mu_i\varepsilon^{\alpha_i}$, где $\varphi_0 \neq 0$, $1<p_0<\alpha_1<\alpha_2<\dots$ --- некоторое собственное значение матрицы $A_{\varepsilon}$ при достаточно малых $\varepsilon$. Согласно теоремы Лапласа имеем
$$
\begin{array}{c}
|A_\varepsilon-\lambda_0^\varepsilon E|=\left| A_\varepsilon-(\varphi_0\varepsilon^{\rho_0}+o(\varepsilon^{\rho_0}))E \right|=\\
=\varepsilon^s(-1)^{p-s}\lambda_{p+1}\lambda_{p+2}\dots\lambda_{n}\begin{vmatrix}
b_{q_11} & b_{q_1p_1} & \cdots & b_{q_1p_{s-1}}\\b_{q_21} & b_{q_2p_1} & \cdots & b_{q_2p_{s-1}}\\ \vdots & \vdots & \ddots & \vdots \\ b_{q_s1} & b_{q_sp_1} & \dots & b_{q_sp_{s-1}}
\end{vmatrix}+o(\varepsilon^s), \quad \varepsilon\to 0,
\end{array}
$$
где $y_i\in C, \xi_i>s$ для $i\ge s+1$. С учётом (\ref{kop:eq24}) отсюда получим, что $\left|A_{\varepsilon}-\lambda_0^{\varepsilon}E\right|\neq 0$, что невозможно, поскольку $\lambda_0^{\varepsilon}$ --- собственное значение матрицы $A_{\varepsilon}$. Итак, $0<p_j\le 1$ для всех $1\le j\le p$.

Пусть $R=\max\limits_{p+1\le j\le n}|\lambda_j|+1$, тогда как множество $R^+$ несчётно, а множество $\left\{2k\pi i\lambda_j^{-1},\quad j=p+1,\dots,n, \quad k\in Z \right\}\cup \left\{2k\pi R^{-1}, \quad k\in Z\right\}$ счётно, то существует такое число $T>0$, что $TR\neq 2k\pi$, $T\lambda_j\neq 2k\pi i$, $\forall k\in Z$, $j=p+1,\dots, n$. Пусть $T_{\varepsilon}=\left[T(2\pi \varepsilon)^{-1}\right]2\pi \varepsilon$, где $[a]$ --- целая часть числа $a$, тогда $T_{\varepsilon}=T+O(\varepsilon)$ при $\varepsilon \to 0$. Поэтому выполняются следующие соотношения
$$  
\begin{array}{c} 
T_{\varepsilon}\lambda_j^{\varepsilon}=\left(T+O(\varepsilon)\right)\left(\varphi_j\varepsilon^{p_j}+o(\varepsilon^{p_j})\right)=\varphi_jT\varepsilon^{p_j}+o(\varepsilon^{p_j}) ,\quad \varepsilon \to 0,\\
j=1,\dots,p,\\
\\T_{\varepsilon}\lambda_j^{\varepsilon}=\left(T+O(\varepsilon)\right)\left(\lambda_j+o(1)\right)=T\lambda_j+o(1), \quad \varepsilon \to 0,\quad j=p+1,\dots,n,\\
\\\\T_{\varepsilon}R=\left(T+O(\varepsilon)\right)R=TR+O(\varepsilon),\quad \varepsilon \to 0,
\end{array}
$$
и при достаточно малых $\varepsilon$ имеют место неравенства $T_{\varepsilon}R\neq 2k\pi$, $ T_{\varepsilon}\lambda_j^{\varepsilon}\neq 2k\pi i$, $\forall k\in Z$, $j=1,\dots,n$. Оценим теперь величину $\parallel \left(e^{T_{\varepsilon}A_{\varepsilon}}-E\right)^{-1}\parallel$.

\begin{lemma}
Для достаточно малых $\varepsilon$ выполняется неравенство $\parallel\left(e^{T_{\varepsilon}A_{\varepsilon}}-E\right)^{-1}\parallel <L\varepsilon^{-1}$ при некотором $L=const>0$.
\end{lemma}

Действительно, пусть $\varphi=\min\limits_{1\le j\le p}|\varphi_j|$, $\Gamma_1=\left\{\lambda\in C: |\lambda |=\frac{1}{2}\varphi \varepsilon \right\}$, $\Gamma_2=\left\{\lambda\in C: |\lambda |=R\right\}$. Тогда при достаточно малых $\varepsilon$ все собственные значения матрицы $A_{\varepsilon}$ лежат внутри области $\Omega$, ограниченной котурами $\Gamma_1$ и $\Gamma_2$. Положим $f(A)=\left(e^{T_{\varepsilon}A}-E\right)^{-1}$ и пусть контуры $\Gamma_1$ и $\Gamma_2$ ориентированы положительно. т.е. так, что при обходе вдоль области остаётся слева. Тогда (см. \cite{kop:bib7}, c. 110)
$$
\begin{array}{c} 
f(A_{\varepsilon})=\frac{1}{2\pi i}\int\limits_{\Gamma_1}\left(e^{T_{\varepsilon}\lambda}-1\right)^{-1}\left(\lambda E-A_{\varepsilon}\right)^{-1}d\lambda +\\
+\frac{1}{2\pi i}\int\limits_{\Gamma_2}\left(e^{T_{\varepsilon}\lambda}-1\right)^{-1}\left(\lambda E-A_{\varepsilon}\right)^{-1}d\lambda 
\end{array}
$$
Заметим здесь, что $\left(\lambda E-A_{\varepsilon}\right)^{-1}=\left(det(\lambda E-A_{\varepsilon})\right)^{-1}(A_{ij})$, где $(A_{ij})$ --- транспонированная матрица алгебраических дополнений элементов матрицы $\lambda E-A_{\varepsilon}$. Для доказательства леммы будем использовать следующие утверждения

\begin{statement}
Для всех $\lambda \in \Gamma_1$ и при достаточно малых $\varepsilon$ существует такое число $L_1>0$, что выполняется неравенство 
\begin{equation}\label{kop:eq25}
\parallel (\lambda E-A_{\varepsilon})^{-1}\parallel \le L_1\varepsilon^{-1}.
\end{equation}
\end{statement}
Доказательство. Пусть $\lambda \in \Gamma_1$, тогда $\lambda=\beta_0\varepsilon$, где $\beta_0 \neq 0$ --- некоторая постоянная. Рассмотрим случай, когда все клетки Жордана $K_r (r=1,\dots, s)$ жордановой формы матрицы $A_0$ являются кратными. т.е. $e_1, e_2,\dots, e_s>1$. Имеем тогда при $\varepsilon\to 0$
$$
det(\lambda E-A_{\varepsilon})=\varepsilon^s(-1)^{n-p+s}\lambda_{p+1}\lambda_{p+2}\dots \lambda_n \begin{vmatrix} b_{{q_1}1} & b_{{q_1}{p_1}} & \cdots & b_{{q_1}{p_{s-1}}}\\
b_{{q_2}1} & b_{{q_2}{p_1}} & \cdots & b_{{q_2}{p_{s-1}}}\\
\vdots & \vdots & \ddots & \vdots \\
b_{{q_s}1} & b_{{q_s}{p_1}} & \cdots & b_{{q_s}{p_{s-1}}} \end{vmatrix} + o(\varepsilon^s).
$$
Положим 
$$
P=\begin{pmatrix}b_{{q_1}1} & b_{{q_1}{p_1}} & \cdots & b_{{q_1}{p_{s-1}}}\\
b_{{q_2}1} & b_{{q_2}{p_1}} & \cdots & b_{{q_2}{p_{s-1}}}\\
\vdots & \vdots & \ddots & \vdots \\
b_{{q_s}1} & b_{{q_s}{p_1}} & \cdots & b_{{q_s}{p_{s-1}}} \end{pmatrix},
$$ 
и обозначим через $P_{ij}(i=q_1,\dots , q_s, j=1,p_1,\dots , p_{s-1})$ --- алгебраические дополнения элементов матрицы $P$. Отсюда
$$
A_{ij}=\varepsilon^{s-1}(-1)^{n-p+s-1}\lambda_{p+1}\lambda_{p+2}\dots \lambda_nP_{ij}+o(\varepsilon^{s-1}),\quad \varepsilon \to 0,
$$
для всех $i=q_1,\dots , q_s, j=1, p_1,\dots , p_{s-1})$, а 
$$
\begin{array}{c} 
A_{ij}=o(\varepsilon^{s-1}),\quad \varepsilon \to 0,\quad (i,j) \notin \left(\{q_1,q_2,\dots,q_s\},\{1,p_1,p_2,\dots,p_{s-1}\}\right).
\end{array}
$$
Так как в силу (\ref{kop:eq24}) $\det P \neq 0$, то существует хотя бы одно алгебраическое дополнение $P_{ij}\neq 0$. Отсюда вытекает неравенство
\begin{equation}\label{kop:eq26}
\begin{array}{c}
\parallel (\lambda E-A_{\varepsilon})^{-1}\parallel \le \frac{1}{\left|\det(\lambda E-A_{\varepsilon})\right|}\left(\sum\limits_{i,j=1}^n|A_{i,j}|^2\right)^{\frac{1}{2}}\le \frac{L_1}{\varepsilon},
\end{array}
\end{equation}
при некотором $L_1=const > 0$ и для достаточно малых $\varepsilon$.

Рассмотрим теперь случай, когда среди клеток Жордана $K_r(r=1,\dots,s)$ жордановой формы матрицы $A_0$ существует некоторые простые клетки $K_{r_1}, K_{r_2}, \dots ,K_{r_l}$ $(1\le r_1<r_2<\dots <r_l\le s$. Положим 
$$
P'=\begin{pmatrix}b_{{q_1}1}-\sigma_1 & b_{{q_1}{p_1}} & \cdots & b_{{q_1}{p_{s-1}}}\\
b_{{q_2}1} & b_{{q_2}{p_1}}-\sigma_2 & \cdots & b_{{q_2}{p_{s-1}}}\\
\vdots & \vdots & \ddots & \vdots \\
b_{{q_s}1} & b_{{q_s}{p_1}} & \cdots & b_{{q_s}{p_{s-1}}}-\sigma_s \end{pmatrix},
$$
где $\sigma_k=\beta_0$, если $k\in\{r_1, r_2, \dots, r_l\}$, а иначе $\sigma_k = 0$, и обозначим через $P'_{ij}$  ($i=q_1,q_2,\dots ,q_s$, $j=1,p_1,p_2,\dots, p_{s-1}$) --- алгебраические дополнения элементов матрицы $P$. Тогда имеют место следующие соотношения
$$
\begin{array}{c}
\det(\lambda E-A_{\varepsilon})=\varepsilon^s(-1)^{n-p+s}\lambda_{p+1}\lambda_{p+2}\dots \lambda_n \det P+o(\varepsilon^s),\quad \varepsilon \to 0,
\end{array}
$$
$$
\begin{array}{c}
A_{ij}=\varepsilon^{s-1}(-1)^{n-p+s-1}\lambda_{p+1}\lambda_{p+2}\dots \lambda_n P'_{ij}+o(\varepsilon^{s-1}), \quad \varepsilon \to 0,
\end{array}
$$
для всех $i=q_1,q_2,\dots ,q_s$, $j=1,p_1,p_2,\dots ,p_{s-1}$, а
$$
\begin{array}{c}
A_{ij}=O(\varepsilon^{s-1}),\quad \varepsilon \to 0,\quad \forall (i,\quad j)\notin \left(\{q_1,q_2,\dots ,q_s\},\quad \{1,p_1,p_2,\dots ,p_{s-1}\}\right).
\end{array}
$$
Поскольку $\det P\neq 0$, а радиус окружности $\Gamma_1$ можно немного повернуть, то, не нарушая общности, можно считать, что $\det P' \neq 0$

т.е. $\sum\limits_{i=0}^{\infty}\beta_i\varepsilon^{i+1}$ является собственным значением матрицы $A_{\varepsilon}$. Для нахождения чисел $\beta_i$, $i=1,2,\dots$, положим $f_i=0, i=1,2,\dots$ Очевидно, что если для некоторого $s\ge 1$ известны все $\beta_i$, $i=0,\dots ,s-1$, то $f_s=0$ есть уравнение первой степени для $\beta_s$. Так как $\beta_0$ --- известное число, то отсюда последовательно решив уравнения $f_s=0$, $s=1,2,\dots$ находим числа $\beta_i$, $i=1,2,\dots$. Следовательно, число $\beta_0$ совпадает с каким-то числом $\varphi_j$, $i=1,\dots,p$, что невозможно в силу того, что $\lambda \in \Gamma_1$. Итак, $\det P' \neq 0$. Отсюда следует неравенство (\ref{kop:eq25}). Утверждение 1 полностью доказано.  
\begin{statement}
Для всех $\lambda \in \Gamma_2$ и при достаточно малых $\varepsilon$ выполняется неравенство $\parallel (\lambda E-A_{\varepsilon})^{-1}\parallel \le L_2$ при некотором $L_2=const>0$.
\end{statement}
Действительно, при $\lambda \in \Gamma_2$ имеем
$$
\begin{array}{c}
det(\lambda E-A_{\varepsilon})=\lambda^p \prod\limits_{j=p+1}^n(\lambda - \lambda_j)+O(\varepsilon) при \varepsilon \to 0,
\end{array}
$$
$$
\begin{array}{c}
A_{ii}=\lambda^{p-1}\prod\limits_{j=p+1}^n(\lambda - \lambda_j)+O(\varepsilon), i=1,\dots,p,
\end{array}
$$
$$
\begin{array}{c} 
A_{ii}=\lambda^p\prod\limits_{j=p+1, j\ne 0}^n(\lambda - \lambda_j)+O(\varepsilon), при \varepsilon \to 0, i=p+1,\dots,n,
\end{array}
$$
$$
\begin{array}{c}
A_{ij}=C_{ij}+O(\varepsilon), \varepsilon \to 0,
\end{array}
$$
где $C_{ij}$ --- некоторые постоянные, $i, j=1,\dots, n$, $i\neq j$.
Отсюда следует
$$
\parallel (\lambda E-A_{\varepsilon})^{-1}\parallel \le \frac{1}{\left|\det(\lambda E-A_{\varepsilon})\right|}\left(\sum\limits_{i,j=1}^n|A_{i,j}|^2\right)^{\frac{1}{2}}\le L_2,\quad L_2=const>0.
$$
Утверждение 2 доказано.

Далее для $\lambda \in \Gamma_1$ по правилу Лопиталя $\frac{e^{|T_{\varepsilon}\lambda |}-1}{|T_{\varepsilon}\lambda|} \to 1$ при $\varepsilon \to 0$. Поэтому при достаточно малых $\varepsilon$ имеет место неравенство $\left|e^{T_{\varepsilon}\lambda }-1)^{-1}\right| \le L_3\varepsilon^{-1}$, где $L_3 > 0$ --- некоторое постоянное. Имеем теперь $T_{\varepsilon}R-2l\pi=(T+O(\varepsilon))R-2l\pi=TR-2l\pi+O(\varepsilon)$ при $\varepsilon \to 0$. Следовательно, для достаточно малых $\varepsilon$, $|T_{\varepsilon}R-2l\pi |>\frac{\delta}{2}>0$, $ \forall l \in Z $, где $\sigma= \min \{TR-2l\pi, l \in Z\}$. Поэтому $|T_{\varepsilon}\lambda-2l\pi |>\frac{\delta}{2}>0$, $ \forall \lambda \in \Gamma_2$, $ \forall l \in Z $. Тогда существует $\sigma = \sigma (\delta)>0$ такое, что $\left|e^{T_{\varepsilon}\lambda}-1\right|>\sigma $ для достаточно малых $\varepsilon$, т.е. существует такое число $L_4>0$, что $\left|(e^{T_{\varepsilon}\lambda}-1)^{-1}\right|\le L_4$, $ \forall \lambda \in \Gamma_2$. В силу утверждений 1 и 2 отсюда находим
$$
\begin{array}{c}
\parallel (e^{T_{\varepsilon}A_{\varepsilon}}-E)^{-1}\parallel = \parallel f(A_{\varepsilon}) \parallel \le\\
\le \parallel \frac{1}{2\pi i}\int\limits_{\Gamma_1}(e^{T_{\varepsilon}\lambda}-1)^{-1}(\lambda E-A_{\varepsilon})^{-1}d\lambda \parallel +\\
+ \parallel \frac{1}{2\pi i}\int\limits_{\Gamma_2}(e^{T_{\varepsilon}\lambda}-1)^{-1}(\lambda E-A_{\varepsilon})^{-1}d\lambda \parallel \le\\
\le \left|(e^{T_{\varepsilon}\lambda}-1)^{-1}\right| \parallel (\lambda E-A_{\varepsilon})^{-1}\parallel \left| \int\limits_{\Gamma_1}d\lambda \right| +\\
+\left|(e^{T_{\varepsilon}\lambda}-1)^{-1}\right| \parallel (\lambda E-A_{\varepsilon})^{-1}\parallel \left| \int\limits_{\Gamma_2}d\lambda \right| \le\\
\le L_1\varepsilon^{-1}L_3\varepsilon^{-1}2\pi \frac{1}{2}|\varphi_k|\varepsilon +L_2L_42\pi R \le \frac{L_5}{\varepsilon}+L_6 \le L\varepsilon^{-1},  
\end{array}
$$  
где $L_5$, $L_6$, $L$ --- некоторые положительные постоянные. Лемма 3 доказана. 

Докажем теперь, что при достаточно больших $\omega$ система (\ref{kop:eq1}) имеет не более одного $T_{\omega^{-1}}$-периодического решения. Пусть $x^{(1)}(t)$, $x^{(2)}(t)$ --- два $T_{\omega^{-1}}$-периодических решения системы (\ref{kop:eq1}), $x^{(0)}(t)=x^{(1)}(t)-x^{(2)}(t)$. Тогда вектор-функция $x^{(0)}(t)$ удовлетворяет уравнению
$$
\frac{dx^{(0)}}{dt}=\left(A_0+\frac{1}{\omega}B_0\right)x^{(0)}+\sum\limits_{1\le |l|\le m}B_lx^{(0)}e^{il\omega t}.
$$  
Производя в последнем последовательно замены переменных $x^{(0)}(t)=y(t) + \sum\limits_{1\le |l| \le m}\frac{B_ly(t)}{il\omega}e^{il\omega t}$, $y(t)=z(t)+\sum\limits_{1\le |l| \le m}\frac{D_lz(t)}{il\omega^2}e^{il\omega t}$, придём (см. \cite{Do1}) к уравнению 
\begin{equation}\label{kop:eq27}
\begin{array}{c}
\frac{dz}{dt}=A_{\omega^{-1}}z+\frac{1}{\omega^2}Q(t, \omega)z,
\end{array}
\end{equation}
где $Q(t,\omega)$ --- $\frac{2\pi}{\omega}$-периодическая относительно $t$ матрица-функция с равномерно ограниченной нормой: $\parallel Q(t,\omega)\parallel \le C_0$, $C_0 > 0$ не зависит от $\omega \gg 1$, $t \in R$. Как известно (см. \cite{kop:bib8}, с. 34), $T_{\omega^{-1}}$-периодическое решение $z(t)$ системы (\ref{kop:eq27}) удовлетворяет уравнению 
$$
\begin{array}{c}
z(t)=\frac{1}{\omega^{k+1}}\left[(E-e^{T_{\omega^{-1}}A_{\omega^{-1}}})^{-1}\int\limits_{0}^{T_{\omega^{-1}}}e^{(T_{\omega^{-1}}-s)A_{\omega^{-1}}}Q(s,\omega)z(s)ds+\right.\\
\left.+\int\limits_{0}^{t}e^{(t-s)A_{\omega^{-1}}}Q(s,\omega)z(s)ds \right].
\end{array}
$$  
Заметим, что при $\omega \ge 1$, $|T_{\omega^{-1}}|=\left|\left[T\frac{\omega}{2\pi}\right]\frac{2\pi}{\omega}\right|\le T$, $\parallel T_{\omega^{-1}}A_{\omega^{-1}}\parallel \le C_1$, где $C_1$ --- некоторое положительное число. Тогда при достаточно больших $\omega$ в силу леммы 3 имеем $\forall t\in [0,T_{\omega^{-1}}]$
$$
\begin{array}{c}
|z|= \parallel z(t) \parallel \le\\
\le \frac{1}{\omega^2}\parallel \left(E-e^{T_{\omega^{-1}}A_{\omega^{-1}}}\right)^{-1}\int\limits_{0}^{T_{\omega^{-1}}}e^{(T_{\omega^{-1}}-s)A_{\omega^{-1}}}Q(s,\omega)z(s)ds\parallel +\\
+ \frac{1}{\omega^2}\parallel \int\limits_{0}^{t}e^{(t-s)A_{\omega^{-1}}}Q(s,\omega)z(s)ds \parallel \le \frac{C_3|z|}{\omega}+\frac{C_4|z|}{\omega^2}\le\frac{C}{\omega}|z|,
\end{array}
$$
где $C_3$, $C_4$, $C$ --- некоторые положительные постоянные, независимые от $\omega$. Отсюда следует, что $z\equiv 0$ при больших $\omega$, т.е. при достаточно больших $\omega$ система (\ref{kop:eq1}) имеет не более одного $T_{\omega^{-1}}$ --- периодического решения.               

\section{Линейная система в частным производных с оператором Стокса в главной части}

\subsection{Основной результат}

Пусть $\Omega$ ~--- ограниченная область в $R^3$ со сколь угодно гладкой границей $\partial\Omega$, $m\in N$, $\omega\gg 1$. В бесконечном цилиндре $Q=\Omega\times R$ рассмотрим задачу о вещественных $\frac{2\pi}{\omega}$-периодических по времени $t$ решениях системы уравнений
\begin{equation}\label{eq:main}
\frac{\partial u}{\partial t} + \nabla p = \Delta u + B_0(x)u + \frac{1}{\omega}C(x)u + \sum\limits_{1\le |k|\le m}\left(L_k(x)u+d_k(x)\right)e^{ik\omega t}+d_0(x),
\end{equation}
\begin{equation}\label{eq:div}
\Div u = 0
\end{equation}
с граничным условием
\begin{equation}\label{eq:bound}
u|_{\partial\Omega} = 0.
\end{equation}
Здесь $x=(x_1,x_2,x_3)\in\Omega$, $t\in R$, $u=u(x,t)$ ~--- неизвестная трехмерная функция, а $B_0(x)$, $C(x)$, $L_k(x)$ и $d_0(x)$, $d_k(x)$ ~--- известные бесконечно гладкие матрицы-функции и вектор-функции соответственно, причем $B_0(x)$, $C(x)$ и $d_0(x)$ ~--- вещественные, а $L_k(x)$, $d_k(x)$ комплексно сопряжены с $L_{-k}(x)$, $d_{-k}(x)$ соответственно.  

Символом $S_2(\Omega)$ будем обозначать замыкание по норме $L_2(\Omega)$ множества непрерывно дифференцируемых комплекснозначных вектор-функций $u$, имеющих на $\partial\Omega$ равную нулю нормальную компоненту и удовлетворяющих условию $\Div u = 0$. Символом $\Pi$ обозначим известный ортогональный проектор в $L_2(\Omega)$ на $S_2(\Omega)$ (см. \cite{Ud1},\cite{Lad1}).  

Введем теперь действующий в $S_2(\Omega)$ оператор $A=\Pi\Delta + \Pi B_0(x)$, область определения $D(A)$ которого является замыканием по норме $W^2_2(\Omega)$ линеала гладких соленоидальных ($div(u)=0$) вектор-функций, обращающихся в нуль на $\partial\Omega$. Определим еще дифференциальное выражение $B=\Pi C(x)+\sum\limits_{1\le |k|\le m}\frac{\Pi L_k(x)\Pi L_{-k}(x)}{ik}$. Предположим, что $\lambda = 0$ ~--- $n$-кратное ($n\in N$) собственное значение оператора $A$, которому отвечает ортонормированный набор из $s (1 \le s \le n)$ собственных функций: $e_1(x)$, $e_2(x)$, ... , $e_s(x)$. Будем предполагать, что ни одна собственная функция $e(x)$ оператора $A$, отвечающая нулевому собственному значению, не имеет присоединенных относительно пары операторов $A$, $B$ \cite{kop:bib4}, т. е. задача 
\begin{equation}
Av(x)=-Be(x),
\end{equation}
\begin{equation}
v(x)|_{\partial\Omega}=0,
\end{equation}
не имеет классических решений, где $e(x)$ ~--- любая собственная функция оператора $A$. Согласно альтернативе Фредгольма, последнее равносильно соотношению
\begin{equation}\label{eq:razr}
\sum\limits_{j=1}^s|(Be(x),b_j(x))|\neq 0,
\end{equation}
где $b_j(x)$, $j=1,2,\dots,s$, ~--- линейно независимые решения уравнения $A^*z(x)=0$. Здесь $A^*$ ~--- сопряженный к $A$, а $(\cdot,\cdot)$ ~--- скалярное произведение в $L_2(\Omega)$. Покажем, что матрица $P\equiv\{(Be_i(x),b_j(x))\}_{i,j=1}^s$ невырождена. Действительно, составим линейную комбиницию из столбцов матрицы $P$ и приравняем ее к нулю:
\begin{equation}\label{eq:stolbcy}
\sum\limits_{i=1}^{s}c_i(Be_i(x),b_j(x))=0,\quad j=1,2,\dots,s.
\end{equation}  
Если набор $c_1$, $c_2$, ... , $c_s$ нетривиален, то из уравнения (\ref{eq:stolbcy}) следует нарушения условия (\ref{eq:razr}), значит линейная комбинация столбцов обращается в ноль только при тривиальном наборе коэффициентов, а матрица $P$ невырождена. 

C целью использования в дальнейшем метода пограничного слоя \cite{Vishik2} перейдем в некоторой окрестности $\Omega_0$ границы $\partial\Omega$ области $\Omega$ к криволинейным координатам следующим образом. Через каждую точку $x\in\partial\Omega$ проведем внутреннюю нормаль, причем окрестность $\Omega_0$ будем считать настолько малой, что нормали в ней не пересекаются. В окрестности $\Omega_0$ положим $x=(\psi_1,\psi_2,r)$, где $r$ расстояние от $x$ до границы $\partial\Omega$ по нормали, $\psi=(\psi_1, \psi_2)$ ~--- координаты соответствующей точки границы.

Асимптотику вещественного $2\pi\omega^{-1}$-периодического решения задачи (\ref{eq:main})-(\ref{eq:bound}) будем искать в виде
\begin{equation}\label{eq:asymp}
\begin{array}{c}
u_\omega(x,t)=\omega \sum\limits_{j=1}^{s}c^{j}_{-2}e_j(x)+\sum\limits_{k =-1}^{\infty}\omega^{-\frac{k}{2}}[u_k(x)+v_k(\psi,\rho)+\\
+\sum\limits_{j=1}^{s}c^{j}_{k}e_j(x)+y_k(x,\omega t)+z_k(\psi,\rho,\omega t)],
\end{array}
\end{equation}
\begin{equation}\label{eq:asympU}
\begin{array}{c}
p_\omega(x)=\sum\limits_{k=-1}^{\infty}\omega^{-\frac{k}{2}}[p_k(x)+\omega^{\frac{1}{2}}s_{k-1}(\psi,\rho)+\omega m_{k-2}(x,\omega t)+\omega^{\frac{1}{2}}n_{k-1}(\psi,\rho,\omega t)].
\end{array}
\end{equation}
Здесь $u_k$, $y_k$, $p_k$, $m_k$ ~--- регулярные, а $v_k$, $z_k$, $s_k$, $n_k$ ~--- погранслойные вектор-функции (см. \cite{Vishik2}), причем $y_k$, $z_k$, $n_k$ являются $2\pi$-периодическими по $\tau$ с нулевым средним, т. е. 
$$
\begin{array}{c}
\langle y_k\rangle(x)\equiv\frac{1}{2\pi}\int\limits_0^{2\pi}y_k(x,\tau)d\tau=0,\\
\langle z_k\rangle(\psi,\rho)=0.
\end{array}
$$
Погранслойные функции вначале определяются при $x\in\Omega_0$, а затем продолжаются нулем вовнутрь области $\Omega$ и умножаются на соответствующую срезающую бесконечно дифференцируемую в $\Omega$ функцию. 

Формулировке основного результата параграфа предпошлем ряд обозначений. 

Задачей (A) назовем задачу Дирихле в $\bar{\Omega}$ вида 
$$
\left\{\begin{array}{c}
(\Delta+ B_0)u(x)+\nabla p(x)=G(x),\\
\Div u(x)=0,\\
(u(x),e_j(x))=0, j=1,2,\dots,s,\\
u(x)|_{\partial\Omega}=0,
\end{array}\right.
$$  
где $G\in C^\infty(\bar{\Omega})$, $(G,b_j)=0$, $j=1,2,\dots,s$. 

Задачей (B) ~--- задачу Неймана в $\bar{\Omega}$
$$
\left\{\begin{array}{c}
\Delta q(x) = g(x),\\
\frac{\partial q(x)}{\partial n}|_{\partial\Omega}=0,
\end{array}\right.
$$
где $\int\limits_{\Omega}g(x)dx=0$, $g\in C^\infty(\bar{\Omega})$. Как известно, задачи (A) и (B) однозначно разрешимы, если единственность $p$ и $q$ понимать с точностью до постоянного слагаемого.

Аналогичным образом, говоря о единственности решения $(u,p)$ системы (\ref{eq:main})-(\ref{eq:bound}), мы подразумеваем, что единственность его компоненты $p$ понимается с точностью до не зависящего от $x$ слагаемого.

Далее, определим частичные суммы $u^n_\omega$ и $p^n_\omega$  рядов (\ref{eq:asymp}), (\ref{eq:asympU}), заменив формально в (\ref{eq:asymp}), (\ref{eq:asympU}) $\infty$ на $n$.

\begin{theorem}
Существует такое число $\omega_0$, что при $\omega>\omega_0$  справедливы следующие утверждения.
1. Задача  (\ref{eq:main})-(\ref{eq:bound})  имеет единственное $\frac{2\pi}{\omega}$-периодическое по $t$ решение $(u_\omega(x,t),p_\omega(x,t))$, которое является вещественным и бесконечно дифференцируемым.
2. Построение вектор-функций $u_\omega^n$, $p_\omega^n$ при каждом $n\ge -1$ сводится к решению конечного (зависящего от $n$) числа задач вида (A) и (B). 
3. Для любого $l\ge 0$ и любого целого $n\ge -1$ справедливы оценки
\begin{equation}\label{eq:ev}
\begin{array}{c}
\|u_\omega - u_\omega^n\|_{C_{x,t}^{l,l/2}}\le c_{n,l}\omega^{-(n+1)+l/2},\\
\|\nabla p_\omega - \nabla p_\omega^n\|_{C_{x,t}^{l,l/2}}\le d_{n,l}\omega^{-n+l/2},
\end{array}
\end{equation}
где $c_{n,l}$, $d_{n,l}=const>0$. 
\end{theorem}

\subsection{Существование, единственность, гладкость и вещественность решения}

Вначале напомним известные вспомогательные результаты (см., например, \cite{Ud1}, \cite{Krasn}, \cite{Simon83}). При этом символом $\|\cdot\|$ будем обозначать норму в пространстве $S_2(\Omega)\equiv H^0$, а символом $Hom(B_1,B_2)$, где  $B_1$, $B_2$ --- банаховы пространства, --- пространство линейных ограниченных операторов, действующих из $B_1$ в $B_2$ с обычной операторной нормой. 

Определенный в п. 2.1 оператор $A$, lдействующий в $H^0$, замкнут и существуют такие положительные числа $a$, $c$, что при $Re\lambda\ge a$ справедлива оценка
$$
\|(\lambda I - A)^{-1}\|\le c(1+|\lambda|)^{-1}. 
$$
Из этой оценки следует, что $A$ порождает в $H^0$ аналитическую полугруппу $e^{tA}\in Hom(H^0, H^0)$, $t\ge 0$, и для любых $T>0$, $\mu\in(0,1)$, найдется постоянная $c$, при которой для всех $u\in H^0$ выполняются неравенства
\begin{equation}\label{eq:oc2}
\left\|\frac{d}{dt}e^{tA}u\right\|\le ct^{-\frac{1}{2}}\|u\|, \quad t\in [0,T], 
\end{equation}
\begin{equation}\label{eq:oc3}
\begin{array}{c}
\|(e^{t_2A}-e^{t_1A})u\|\le ct_1^{-\mu}(t_2-t_1)^\mu \|u\|,\quad 0<t_1<t_2\le T.
\end{array}
\end{equation} 

Подействуем на уравнение (\ref{eq:main}) проектором $\Pi$ и перепишем (\ref{eq:main})-(\ref{eq:bound}) в операторной форме
\begin{equation}\label{eq:oper}
\frac{\partial u}{\partial t} = A u + \frac{1}{\omega}K u + \sum\limits_{1\le |k|\le m}\left(M_k u+a_k(x)\right)e^{ik\omega t}+a_0(x).
\end{equation}
Здесь $A$, $K=\Pi C(x)$ и $M_k=\Pi L_k(x)$ ~--- линейные операторы в $H^0$ с областями определения  
$$
D(A) = \left\{u\in S_2(\Omega)\cap W_2^2(\Omega),\quad u|_{\partial\Omega}=0\right\},
$$
$D(K)=D(M_k)=D(A)$; $a_s(x)=\Pi d_s(x)$ --- вектор-функции со значениями в $H^0$.

Утверждение теоремы 3 о существовании и единственности решения задачи (\ref{eq:main})-(\ref{eq:bound}) эквивалентно аналогичному утверждению для уравнения (\ref{eq:oper}). C доказательства последнего мы и начнем.

 В уравнении  (\ref{eq:oper}) сделаем замену переменных типа замены Крылова-Боголюбова \cite{kop:bib1}
 $$
 u=v+\sum\limits_{0<|k|\le m}(ik\omega I-A-\omega^{-1}K)^{-1}M_ke^{ik\omega t}v\equiv v+S_\omega v.
 $$
 Получим
 $$
 \begin{array}{c}
 (I+S_\omega)\frac{dv}{dt}+\sum\limits_{0<|k|\le m}ik\omega(ik\omega I - A - \omega^{-1} K)^{-1}M_ke^{ik\omega t}v=\\
=Av + \sum\limits_{0<|k|\le m}A(ik\omega I-A-\omega^{-1}K)^{-1}M_ke^{ik\omega t}v +\\
+ \omega^{-1}K v + \omega^{-1}\sum\limits_{0<|k|\le m}K(ik\omega I-A-\omega^{-1}K)^{-1}M_ke^{ik\omega t}v +\\
+ \sum\limits_{0< |k|\le m}M_k e^{ik\omega t}v+ \sum\limits_{0\le |k|,|j|\le m}N_j(ik\omega I-A-\omega^{-1}K)^{-1}M_ke^{i(k+j)\omega t}v +\\
+ \sum\limits_{0\le|k|\le m}a_k(x))e^{ik\omega t}.
 \end{array}
 $$
 Отсюда находим
 \begin{equation}\label{eq:operform1}
 \begin{array}{c}
 \frac{dv}{dt}=(I+S_\omega)^{-1}\left[A+\omega^{-1}K+ \sum\limits_{0< |k|\le m}M_{-k}(ik\omega I-A-\omega^{-1}K)^{-1}M_k\right]v +\\
 +(I+S_\omega)^{-1}\left[\sum\limits_{0< |k|,|j|\le m,k+j\neq 0}N_j(ik\omega I-A-\omega^{-1}K)^{-1}M_ke^{i(k+j)\omega t}\right]v+\\
  + (I+S_\omega)^{-1}\sum\limits_{0\le|k|\le m}a_k(x))e^{ik\omega t}.
 \end{array}
 \end{equation}
 Перепишем (\ref{eq:operform1}) в виде
 \begin{equation}\label{eq:operform2}
  \frac{dv}{dt} - A_\omega v = R_\omega v + \sum\limits_{0\le|k|\le m} (I+S_\omega)^{-1}a_k(x)e^{ik\omega t},
 \end{equation}
где
$$
A_\omega = A+\omega^{-1}\left(K+\sum\limits_{0< |k|\le m}\frac{1}{ik}M_{-k}M_k\right)\equiv A+\omega^{-1}B,
$$
$$
\begin{array}{c}
R_\omega = (I+S_\omega)^{-1}S_\omega\left[A+\omega^{-1}K+ \sum\limits_{0< |k|\le m}M_{-k}(ik\omega I-A-\omega^{-1}K)^{-1}M_k\right]+\\
+(I+S_\omega)^{-1}\left[\sum\limits_{0< |k|,|j|\le m,k+j\neq 0}N_j(ik\omega I-A-\omega^{-1}K)^{-1}M_ke^{i(k+j)\omega t}\right]+\\
+\sum\limits_{0< |k|\le m}M_{-k}[(ik\omega I-A-\omega^{-1}K)^{-1}-(ik\omega)^{-1}]M_k.
\end{array}
$$

Перейдем к исследованию некоторых спектральных характеристик оператора $A_\omega$ при $\omega\gg 1$.

Пусть $r_1$ --- столь малое положительное число, что внутри окружности $\Gamma_1=\{\lambda\in C:|\lambda|=r_1\}$ нет ненулевых собственных значений оператора $A$. Введем спектральные проекторы:
$$
P=\int\limits_{|\lambda|=r_1}(\lambda I-A)^{-1}d\lambda,\quad Q = I-P.
$$
Обозначим через $n$ размерность образа $P$ и отметим, что далее в этом параграфе мы будем существенно использовать матричные представления операторов, связанные с определенными парами взаимно ортогональных проекторов. Представления такого рода активно использовал Л. И. Сазонов в работах \cite{Sazon}, \cite{Sazon2}, \cite{Sazon3}.

Воспользуемся матричным представлением оператора

$$
\begin{array}{c}
A_\omega - \lambda I = \begin{pmatrix}
P(A_\omega-\lambda I)P & P A_\omega Q \\
Q A_\omega P & Q(A_\omega-\lambda I)Q
\end{pmatrix}: \\
D(A)\subset PH^0\times QH^0\to PH^0\times QH^0, \quad \lambda \in C.
\end{array}
$$
Учитывая, что оператор $Q(A_\omega-\lambda I)Q$ при больших $\omega$ и $|\lambda|<r_1$ обратим в $QH^0$, введем в рассмотрение оператор
$$
J = \begin{pmatrix}
P & -C_{0\omega} \\
0 & Q
\end{pmatrix},
$$
где
\begin{equation}\label{eq:matrixform1}
C_{0\omega} = \omega^{-1}PBQ(A_\omega-\lambda I)^{-1}Q.
\end{equation}
При этом
\begin{equation}\label{eq:matrixform2}
\begin{array}{c}
J(A_\omega-\lambda I) =\\
= \begin{pmatrix}
P(A_\omega-\lambda I)P-\omega^{-2}PBQ(A_\omega-\lambda I)^{-1}Q B P & 0 \\
\omega^{-1}Q B P & Q(A_\omega-\lambda I)Q
\end{pmatrix}\equiv\\
\equiv\begin{pmatrix}
C_{1\omega} & 0 \\
C_{2\omega} & C_{3\omega}
\end{pmatrix}\equiv C_{\omega}:D(A)\subset H^0\to H^0
\end{array}
\end{equation}
Далее учтем два следующих простых соображения. При больших $\omega$ и $|\lambda |<r_1$ оператор $J:H^0\to H^0$ обратим, причем 
\begin{equation}\label{eq:obrJ}
J^{-1} = \begin{pmatrix}
P & C_{0\omega} \\
0 & Q
\end{pmatrix},
\end{equation}
а оператор $C_\omega$ обратим тогда и только тогда, когда обратим оператор $C_{1\omega}$, причем в этом случае
\begin{equation}\label{eq:obrComega}
C_\omega^{-1} = \begin{pmatrix}
C_{1\omega}^{-1} & 0 \\
-C_{3\omega}^{-1}C_{2\omega}C_{1\omega}^{-1} & C_{3\omega}^{-1}
\end{pmatrix},
\end{equation}
откуда в силу равенства (\ref{eq:matrixform2}) следует что число $\lambda$ при $|\lambda |<r_1$ является собственным значением оператора $A_\omega$ при больших $\omega$ тогда и только тогда, когда оно является собственным значением конечномерного оператора $PAP+\omega^{-1}P[B-\omega^{-1}BQ(A_\omega-\lambda I)^{-1}QB]P$ в пространстве $PH^0$.

Теперь заметим, что уравнение $Ax=Ba$, $x\in H^0$, $a\in Ker A$ эквивалентно системе
$$
\left\{\begin{array}{l}
PAPx_1=PBPa\\
QAQx_2=QBQa,\quad x_1\in PH^0,\quad x_2 \in QH^0, \quad x=(x_1,x_2)^T.
\end{array}\right.
$$
Поскольку второе уравнение этой системы однозначно разрешимо, то вектор $a\in\Ker A$ не имеет присоединенных относительно пары операторов $A$, $B$ в $H^0$ тогда и только тогда, когда он не имеет присоединенных относительно пары операторов $PAP$, $PBP$ в $PH^0$. Таким образом, вопрос о малых  собственных значениях $\lambda_{\omega j}$ оператора $A_\omega$, $\omega\gg 1$, сведен к рассмотреному в \S 1 конечномерному  случаю, поэтому: 
$$
\lambda_{\omega j}=\varphi_j\omega^{-p_j}+o\left(\omega^{-p_j}\right), \quad j=\overline{1,s},\quad \omega\gg 1, \quad\varphi_j\neq 0,\quad 0<p_j\le 1.
$$

Как отмечалось выше, оператор $A$ порождает аналитическую полугруппу $e^{tA}$, $t\ge 0$, а потому в комплексной плоскости $C$ найдется острый угол $\{\lambda\in C:\arg(\lambda-a_0)=\pm\theta_0 \}$, где $a_0>0$, $\theta_0\in\left(\frac{\pi}{2},\pi\right)$ внутри которого лежит спектр оператора $A$. Пусть $t_0$ --- столь малое положительное число, что при всех ненулевых $\lambda\in C$, лежащих внутри угла $\Gamma_0$, $e^{\lambda t_0}\neq 1$. Определенное выше число  $\omega_0$ будем, по-прежнему, считать достаточно большим. Введем в рассмотрение при $\omega > \omega_0$ окружность $\gamma_\omega=\{\lambda\in C:|\lambda|=r_0\omega^{-1}\}$ и число $T_\omega=\frac{2\pi}{\omega}\left[t_0\frac{\omega}{2\pi}\right]$. Тогда при больших $\omega$ действующий в $H^0$ оператор $I-e^{T_\omega A_\omega}$ обратим и справедливо следующее представление, которое понадобится нам в дальнейшем:
\begin{equation}\label{eq:kontursform}
\begin{array}{c}
Z_\omega(t)=(I-e^{T_\omega A_\omega})^{-1}e^{(T_\omega+t-\tau)A_\omega}=\\
=-\frac{1}{2\pi i}\int\limits_{\Gamma_0}\left(1-e^{\lambda T_\omega}\right)^{-1}e^{\lambda (T_\omega + t - \tau)}(A_\omega-\lambda I)^{-1}d\lambda-\\
-\frac{1}{2\pi i}\int\limits_{\gamma_\omega}\left(1-e^{\lambda T_\omega}\right)^{-1}e^{\lambda (T_\omega + t - \tau)}(A_\omega-\lambda I)^{-1}d\lambda,
\end{array}
\end{equation}
где контуры $\Gamma_0$, $\gamma_\omega$ ориентированы так, что при движении вдоль них спектр оператора $A_\omega$ лежит слева.

Задача о классических $2\pi\omega^{-1}$-периодических по $t$ решениях уравнения (\ref{eq:operform2}) эквивалентна задаче о $T_\omega$-периодических решениях этого же уравнения. Последняя же эквивалентна задаче о решениях $z\in C_\mu([0,t_0],H^0)\equiv C_\mu(H^0)$, $\mu\in(0,1)$, уравнения
\begin{equation}\label{eq:intform}
\begin{array}{c}
z(t)=\int\limits_0^t e^{(t-\tau)A_\omega}\left[R_\omega z(\tau)+\sum\limits_{0\le|k|\le m} (I+S_\omega)^{-1}a_k e^{ik\omega t}\right]d\tau+\\
+\left(I-e^{T_\omega A_\omega}\right)^{-1}\int\limits_0^{T_\omega} e^{(t+T_\omega-\tau)A_\omega}\left[R_\omega z(\tau)+\sum\limits_{0\le|k|\le m} (I+S_\omega)^{-1}a_k e^{ik\omega t}\right]d\tau =\\
=\int\limits_0^t e^{(t-\tau)A_\omega}R_\omega z(\tau)d\tau+\left(I-e^{T_\omega A_\omega}\right)^{-1}\int\limits_0^{T_\omega} e^{(t+T_\omega-\tau)A_\omega}R_\omega z(\tau)d\tau+\Psi_\omega(t)\equiv\\
\equiv K_{\omega 1}z + K_{\omega 2} z + \Psi_\omega \equiv K_{\omega}z + \Psi_\omega, \quad t\in[0,t_0].
\end{array}
\end{equation}

Рассмотрим оператор $K_\omega:C_\mu(H^0)\to C_\mu(H^0)$ при $\omega>\omega_0$. Существование и единственность решения уравнения (\ref{eq:intform}) вытекает из следующего утверждения.

\begin{lemma}
Оператор $K_\omega:C_\mu(H^0)\to C_\mu(H^0)$ при $\omega>\omega_0$, где $\omega_0$ достаточно велико, удовлетворяет оценке
\begin{equation}\label{eq:norm}
\|K_\omega\|_{\Hom(C_\mu(H^0),C_\mu(H^0))}<\frac{1}{2}.
\end{equation}
\end{lemma}

Для доказательства последнего неравенства достаточно установить оценки
\begin{equation}\label{eq:key0}
\|K_{\omega i}\|_{\Hom(C_\mu(H^0),C_\mu(H^0))}<\frac{1}{4},\quad i=1,2.
\end{equation}
Интегральные операторы $K_{\omega i}$, $i=1,2$, оцениваются аналогично, поэтому ограничимся лишь оценкой $K_{\omega 2}$, в котором дополнительно (по сравнению с $K_{\omega 1}$) cодержится операторный сомножитель $\left(I-e^{T_\omega A_\omega}\right)^{-1}$. Начнем с матричного представления резольвенты (см. (\ref{eq:matrixform2})-(\ref{eq:obrComega})):
$$
(A_\omega-\lambda I)^{-1} = \begin{pmatrix}
C_{1\omega}^{-1} & 0 \\
-C_{3\omega}^{-1}C_{2\omega}C_{1\omega}^{-1} & C_{3\omega}^{-1}
\end{pmatrix}\begin{pmatrix}
P & -C_{0\omega} \\ 
0 & Q
\end{pmatrix}.
$$
Для $\forall x=\begin{pmatrix}x_1\\x_2\end{pmatrix}\in PH^0\times QH^0$ находим:
\begin{equation}\label{eq:resolventvector}
\begin{array}{c}
(A_\omega-\lambda I)^{-1}x=\begin{pmatrix}
C_{1\omega}^{-1} & 0 \\
-C_{3\omega}^{-1}C_{2\omega}C_{1\omega}^{-1} & C_{3\omega}^{-1}
\end{pmatrix}\begin{pmatrix}
x_1-C_{0\omega}x_2 \\ 
x_2
\end{pmatrix}=\\
=\begin{pmatrix}
C_{1\omega}^{-1}(x_1-C_{0\omega}x_2) \\
-C_{3\omega}^{-1}C_{2\omega}C_{1\omega}^{-1}(x_1-C_{0\omega}x_2) + C_{3\omega}^{-1}x_2
\end{pmatrix}.
\end{array}
\end{equation}
Из соотношений (\ref{eq:kontursform}), (\ref{eq:resolventvector}) следуют представления:
\begin{equation}\label{eq:key1}
\begin{array}{c}
Z_\omega(t)Px=(I-e^{T_\omega A_\omega})^{-1}e^{(T_\omega+t-\tau)A_\omega}=\\
=-\frac{1}{2\pi i}\int\limits_{\Gamma_0}\left(1-e^{\lambda T_\omega}\right)^{-1}e^{(\lambda T_\omega + t - \tau)}\begin{pmatrix}
C_{1\omega}^{-1}x_1 \\
-C_{3\omega}^{-1}C_{2\omega}C_{1\omega}^{-1}x_1
\end{pmatrix}d\lambda-\\
-\frac{1}{2\pi i}\int\limits_{\gamma_\omega}\left(1-e^{\lambda T_\omega}\right)^{-1}e^{\lambda (T_\omega + t - \tau)}\begin{pmatrix}
C_{1\omega}^{-1}x_1 \\
-C_{3\omega}^{-1}C_{2\omega}C_{1\omega}^{-1}x_1
\end{pmatrix}d\lambda\equiv Z_\omega^{(1)}(t)x+Z_\omega^{(2)}(t)x,
\end{array}
\end{equation}
\begin{equation}\label{eq:key2}
\begin{array}{c}
Z_\omega(t)Qx=(I-e^{T_\omega A_\omega})^{-1}e^{(T_\omega+t-\tau)A_\omega}=\\
=-\frac{1}{2\pi i}\int\limits_{\Gamma_0}\left(1-e^{\lambda T_\omega}\right)^{-1}e^{\lambda (T_\omega + t - \tau)}\begin{pmatrix}
-C_{1\omega}^{-1}C_{0\omega}x_2 \\
(C_{3\omega}^{-1}C_{2\omega}C_{1\omega}^{-1}C_{0\omega} + C_{3\omega}^{-1})x_2
\end{pmatrix}d\lambda-\\
-\frac{1}{2\pi i}\int\limits_{\gamma_\omega}\left(1-e^{\lambda T_\omega}\right)^{-1}e^{\lambda (T_\omega + t - \tau)}\begin{pmatrix}
-C_{1\omega}^{-1}C_{0\omega}x_2 \\
(C_{3\omega}^{-1}C_{2\omega}C_{1\omega}^{-1}C_{0\omega} + C_{3\omega}^{-1})x_2
\end{pmatrix}d\lambda\equiv\\
\equiv Z_\omega^{(3)}(t)x+Z_\omega^{(4)}(t)x,\quad t\in [0, t_0
],\quad \tau\in [0,T_\omega],
\end{array}
\end{equation}
Оценка (\ref{eq:key0}) для оператора $K_{\omega 2}$, определяемого равенством (\ref{eq:intform}), осуществляется по хорошо разработанной в теории метода усреднения схеме (см., например, \cite{Simon83}, \cite{Lev7} и формулы (\ref{eq:oc2}), (\ref{eq:oc3})). При этом стандартным образом оценивается оператор $R_\omega$ и используются следующие две леммы.

\begin{lemma}
При $\mu\in(0,1)$ для любых $t\in [0, t_0]$, $\tau\in [0,T_\omega]$ и $\omega>\omega_0$, где $\omega_0$ --- достаточно большое число, справедливы оценки: 
$$
\|Z_\omega^{(i)}(t)x\|_{H^0}\le c\|x\|_{H^0}, \quad \left\|\frac{d}{dt}Z_\omega^{(k)}(t)x\right\|_{H^0}\le c\|x\|_{H^0},\quad i=1,3,4,\quad k=1,4;
$$
$$
\left\|Z_\omega^{(3)}(t_2)x-Z_\omega^{(3)}(t_1)x\right\|_{H^0}\le c(T_\omega + t_1 -\tau)^{-\mu}(t_2-t_1)^\mu, 0<t_1<t_2\le T_0;
$$
$$
\|Z_\omega^{(2)}(t)x\|_{H^0}\le c\omega\|x\|_{H^0}, \quad \left\|\frac{d}{dt}Z_\omega^{(2)}(t)x\right\|_{H^0}\le c\omega\|x\|_{H^0},
$$
где $c$ --- не зависящая от $t$, $\omega$ и $x$ постоянная.
\end{lemma}

\begin{lemma}
При любом $t\in [t,t_0]$ и $\omega>\omega_0$, где $\omega_0$ --- достаточно большое число, справедлива оценка: 
\begin{equation}\label{eq:estimate1}
\|S_\omega A_0\|_{\Hom(H^0, H^0)}\le c_0<\infty,
\end{equation}
где $c_0$ --- не зависящая от $t$, $\omega$ постоянная. 
\end{lemma}

При доказательстве леммы 5 используются представления (\ref{eq:key1}), (\ref{eq:key2}) операторов $Z_\omega^{(i)}$, $i=\overline{1,4}$ и (\ref{eq:matrixform1}), (\ref{eq:matrixform2}) операторов $C_{k\omega}$, $k=\overline{0,3}$. Кроме того, при оценке конечномерного оператора $Z_{\omega}^{(2)}(t)$ используется та же, что в п. 1.4, методика.

Для доказательства леммы 6, очевидно, достаточно установить оценку:
\begin{equation}\label{eq:estimate2}
\|(ik\omega I-A_0)^{-1}M_k A_0\|_{\Hom(H^0, H^0)}\le c_1,\footnote{Оцениваемый в (\ref{eq:estimate2}) оператор является продолжением на $H^0$ одноименного оператора, определенного на линеале $D(A_0)\subset H^0$}
\end{equation}
где $A_0=\Pi\Delta$, а $c_1$ --- независящая от $t$, $\omega$ и $k$, $0<|k|\le m$ постоянная. С этой целью введем в рассмотрение множество $M_0$ гладких соленоидальных финитных в $\Omega$ вектор-функций. Согласно лемме 5.2 \cite{Ud1} множество $M_0$ плотно в $S_2$. При этом для любых $u,v\in M_0$ имеем:
$$
\begin{array}{c}
\left((ik\omega I-A_0)^{-1}M_k A_0u,v\right)=-\left(\Pi\Delta u,M_k^*(ik\omega I+A_0)^{-1}v\right)\equiv\\
\equiv (\Delta u, c_{k\omega}+\grad q_{k\omega}),
\end{array}
$$ 
где 
$$
\Delta q_{k,\omega}=-\Div c_{k\omega}v,\quad \left.\frac{\partial q_{k\omega}}{\partial n}\right|_{\partial\Omega}=0.
$$
Откуда следует, что
$$
\left|\left(ik\omega I-A_0)^{-1}M_k A_0u,v\right)\right|=|(u, c_{k\omega}+\Delta(\grad q_{k\omega}))|\le c_1\|u\|_{H^0}\|v\|_{H^0},
$$ 
где $c_1=const$, таким образом неравенства (\ref{eq:estimate1}), (\ref{eq:estimate2}) доказаны.

Из неравенства (\ref{eq:norm}), в силу принципа сжатых отображений, следует существование единственного решения $z\in C_{\mu}([0,T_\omega],H^0)$ уравнения (\ref{eq:intform}) при больших $\omega$, т. е. $z(t)$ является слабым решением задачи (\ref{eq:main})-(\ref{eq:bound}) с начальным условием $z|_{t=0}=z(0)$. Слагаемые уравнения (\ref{eq:intform}), содержащие $\frac{d}{dt}$ проинтегрируем по частям. После этого выберем и подставим в правую часть последовательность $\{z^{(n)}\}_{n=1}^{\infty}\in C_{\mu}([0,T_\omega],H^0\cap W_{2,0}^2(\Omega))$, которая сходится к $z$ по норме $H^0$. Здесь и далее символом $W_{p,0}^{j}(\Omega)$ обозначается замыкание по норме $W_{p}^{j}(\Omega)$ гладких вектор-функций, обращающихся в нуль на границе $\partial\Omega$. Тогда для соответствующей последовательности $\{\bar{z}^{(n)}\}_{n=1}^{\infty}$ левых частей получим
$$
\bar{z}^{(n)}\in C_{\mu}([0,T_\omega],H^0(\Omega)\cap W_{2,0}^2(\Omega)),\quad \|\bar{z}^{(n)} - z\|_{C_{\mu}([0,T_\omega],H^0(\Omega)\cap W_{2}^1(\Omega))}\to 0.
$$
Таким образом, $z\in C_{\mu,0}([0,T_\omega],H^0(\Omega)\cap W_{2,0}^1(\Omega))$. Из \cite{Solon} теперь следует, что 
\begin{equation}\label{eq:eq14}
z\in C_{\mu}([0,T_\omega],H^0\cap W_{2}^2(\Omega))\cap C_{\mu}([0,T_\omega],H^0\cap W_{2,0}^1(\Omega)),
\end{equation}
поэтому $z(t)$ удовлетворяет задаче (\ref{eq:main})-(\ref{eq:bound}) в $C_{\mu}([0,T_\omega],H^0)$. 

Обозначим символом $C_\mu(R,B)$ обычное банахово пространство вектор-функций $z:R\to B$, удовлетворяющих равномерному условию Гельдера с показателем $\mu$. Продолжим вектор-функцию $z(t)$ $T_\omega$-периодическим образом на всю ось $t\in R$, оставим за продолжением прежнее обозначение. Из (\ref{eq:eq14}) следует, что
\begin{equation}\label{eq:eq15}
z\in C_{\mu}(R,H^0\cap W_{2}^2(\Omega))\cap C_{\mu}(R,H^0\cap W_{2,0}^1(\Omega)).
\end{equation}
В силу $\frac{2\pi}{\omega}$-периодичности коэффициентов исходной задачи (\ref{eq:main})-(\ref{eq:bound}) и доказанной выше единственности решения задачи (\ref{eq:main})-(\ref{eq:bound}) в смысле $C_{\mu}([0,T_\omega],H^0)$, $z(t)$ является единственным $\frac{2\pi}{\omega}$-периодическим решением указанной задачи. Из соотношения (\ref{eq:norm}) следует оценка при больших $\omega$:
\begin{equation}\label{eq:ost}
\|z\|_{C_{\mu}(R,H^0\cap W_{2,0}^1(\Omega))}\le c\sum\limits_{0\le k\le m}\|a_k\|_{H^0}.
\end{equation} 

Для завершения доказательства утверждения 1 осталось установить бесконечную дифференцируемость решения $z(x,t)$: $z\in C^\infty(\bar{Q})$. Пусть $\chi(t)$, $t\in[0,1]$, ~---  бесконечно-дифференцируемая функция такая что $\chi(t)\neq 0$, $t\in(\frac{1}{3},\frac{2}{3})$ и
$$
\chi(t)=\left\{\begin{array}{l}
0,\quad t\in[0,\frac{1}{3}],\\
1,\quad t\in[\frac{2}{3},1].\\
\end{array}\right.
$$
Вектор-функция $v(t)=\chi(t)z(t)$, $t\in[0,1]$, очевидно, является решением задачи Коши
\begin{equation}\label{eq:operdiff}
\begin{array}{c}
\frac{\partial v}{\partial t} - A_\omega v = \chi\left(R_\omega z+\sum\limits_{0\le |k|\le m}(I+S_\omega)^{-1}a_ke^{ik\omega t}\right)+\frac{d\chi}{dt}z,\\
v(0)=0.
\end{array}
\end{equation}
Поскольку коэффициенты левой части (\ref{eq:operdiff}) и граница $\partial\Omega$ бесконечно дифференцируемы, а также для этой задачи выполнены условия согласования сколь угодно высокого порядка, то будем к ней применять теорему 2 из \cite{Solon}. В силу $\frac{2\pi}{\omega}$-периодичности вектор-функции $z(x,t)$ по $t$ получим тогда, что $z\in C^\infty(\overline{Q})$.

\subsection{Построение асимптотики}

Проведем теперь замену $\rho = r\sqrt{\omega}$ и воспользуемся известными (см., например, \cite{Kibel}, \cite{Lev1}) соотношениями для операторов $\Delta$, $\nabla$, $\Div$ в криволинейных координатах $(\psi_1,\psi_2,r)$:
\begin{equation}\label{delta}
(\Delta+B_0) u=\sum\limits_{k=-2}^\infty\omega^{-\frac{k}{2}}N_k(\psi,\rho)u,
\end{equation}
(символ $\Delta$ здесь отождествляется с матричным выражением $\Delta E$, где $E$ ~--- единичная матрица третьего порядка),
\begin{equation}\label{div}
\begin{array}{c}
\Div u=\sum\limits_{k=-1}^\infty\omega^{-\frac{k}{2}}\frac{\partial u^{(3)}}{\partial\rho}-\sum\limits_{k=0}^\infty\omega^{-\frac{k}{2}}[D_{k,1}(\psi,\rho)u^{(1)}+D_{k,2}(\psi,\rho)u^{(2)}\\
+D_{k,3}(\psi,\rho)u^{(3)}],
\end{array}
\end{equation}
\begin{equation}\label{nabla}
\nabla p=\sum\limits_{k=-1}^\infty\omega^{-\frac{k}{2}}P_k(\psi,\rho)p.
\end{equation}
Здесь $P_{-1}(\psi,\rho)p=\left(0,0,\frac{\partial p}{\partial\rho}\right)^T$, $N_{-2}(\psi,\rho)u=\left(\frac{\partial^2 u^{(1)}}{\partial\rho^2},\frac{\partial^2 u^{(2)}}{\partial\rho^2},0\right)^T$, а остальные $P_k$, $N_k$, $D_{j,k}$ ~--- дифференциальные выражения относительно $(\psi,\rho)$ с равномерно ограниченными относительно $(\psi,\rho)$ коэффициентами, $u^{(i)}$, $i=1,2,3$, ~--- $i$-ая криволинейная координата вектора $u$.

Подставим ряды (\ref{eq:asymp}) и (\ref{eq:asympU}) в уравнения (\ref{eq:main}), (\ref{eq:div}) и граничное условие (\ref{eq:bound}), сделав предварительно замену переменной $\tau = \omega t$:
\begin{equation}\label{eq:mainC}
\begin{array}{c}
\sum\limits_{k=-1}^\infty\omega^\frac{2-k}{2}\left(\frac{\partial y_k(x,\tau)}{\partial\tau}+\frac{\partial z_k(\psi,\rho,\tau)}{\partial\tau}\right)+\\
+\sum\limits_{k=-1}^{\infty}\omega^{-\frac{k}{2}}\nabla[p_k(x)+\omega^{\frac{1}{2}}s_{k-1}(\psi,\rho)+\omega m_{k-2}(x,\omega t)+\omega^{\frac{1}{2}}n_{k-1}(\psi,\rho,\omega t)]=\\
=\sum\limits_{k=-1}^\infty\omega^{-\frac{k}{2}}(\Delta+ B_0(x))(u_k(x)+v_k(\psi,\rho)+y_k(x,\tau)+z_k(\psi,\rho,\tau))+\\
+\sum\limits_{j=1}^{s}c^{j}_{-2}C_0(x)e_j(x)+\\
+\sum\limits_{k=-1}^\infty\omega^{-\frac{k+2}{2}} C_0(x)(u_k(x)+v_k(\psi,\rho)+\sum\limits_{j=1}^{s}c^{j}_{k}e_j(x)+y_k(x,\tau)+z_k(\psi,\rho,\tau))+\\
+\sum\limits_{1\le |k| \le m}\left(\omega \sum\limits_{j=1}^{s}c^{j}_{-2}L_k e_j(x)+\sum\limits_{i=-1}^\infty \omega^{-\frac{i}{2}} L_k[u_i(x)+v_i(\psi,\rho)\right.+\\
\left.+\sum\limits_{j=1}^{s}c^{j}_{i}e_j(x)+y_i(x,\tau)+z_i(\psi,\rho,\tau)]+ d_k(x)\right)e^{ik\tau} + d_0(x),
\end{array}
\end{equation}
\begin{equation}\label{eq:divC}
\begin{array}{c}
\sum\limits_{k=-1}^{\infty}\omega^{-\frac{k}{2}}\Div[u_k(x)+v_k(\psi,\rho)+\\
+\sum\limits_{j=1}^{s}c^{j}_{k}e_j(x)+y_k(x,\tau)+z_k(\psi,\rho,\tau)]=0,
\end{array}
\end{equation}
\begin{equation}\label{eq:boundC}
\begin{array}{c}
\sum\limits_{k=-1}^{\infty}\omega^{-\frac{k}{2}}[u_k(x)+v_k(\psi,\rho)+\\
+\sum\limits_{j=1}^{s}c^{j}_{k}e_j(x)+y_k(x,\tau)+z_k(\psi,\rho,\tau)]|_{\partial\Omega}=0.
\end{array}
\end{equation}

В равенствах (\ref{eq:mainC})-(\ref{eq:boundC}) с учетом представлений (\ref{delta})-(\ref{nabla}) приравняем коэффициенты при одинаковых степенях $\omega$ отдельно для регулярных и погранслойных функций. Получим бесконечную последовательность задач. Выпишем первые из них. Из уравнений $(\ref{div})$, $(\ref{eq:divC})$ для погранслойных функций получим
$$
\left\{\begin{array}{c}
\frac{\partial z_{-1}^{(3)}(\psi,\rho,\tau)}{\partial\rho}=0,\\
z_{-1}^{(3)}|_{\rho\to\infty}=0,
\end{array}\right.
$$
$$
\left\{\begin{array}{c}
\frac{\partial v_{-1}^{(3)}(\psi,\rho)}{\partial\rho}=0,\\
v_{-1}^{(3)}|_{\rho\to\infty}=0,
\end{array}\right.
$$
так что $z_{-1}^{(3)}(\psi,\rho,\tau)\equiv 0$, $v_{-1}^{(3)}(\psi,\rho)\equiv 0$. Из уравнений (\ref{eq:mainC})-(\ref{eq:divC}) для регулярных функций находим
\begin{equation}\label{equ1}
\left\{\begin{array}{c}
\frac{\partial y_{-1}(x,\tau)}{\partial\tau}+\nabla m_{-3}(x,\tau)=0,\\
\Div y_{-1}=0,\\
y_{-1n}|_{\partial\Omega}=-z^{(3)}_{-1}|_{\partial\Omega}|_{\rho=0}=0,\\
\langle y_{-1}\rangle=\langle m_{-3}\rangle=0.
\end{array}\right.
\end{equation}
Нижним индексом $n$ мы обозначаем проекцию вектора на внутреннюю нормаль к $\partial\Omega$. Положим $y_{-1}=0$, $m_{-3}=0$. Для погранслойных функций имеем:
$$
\begin{array}{c}
\frac{\partial z_{-1}(\psi,\rho,\tau)}{\partial\tau}=N_{-2}(\psi,\rho)(z_{-1}(\psi,\rho,\tau)+v_{-1}(\psi,\rho))-P_{-1}(s_{-2}(\psi,\rho)+n_{-2}(\psi,\rho,\tau)),\\
\langle z_{-1}\rangle = \langle n_{-2}\rangle = 0.
\end{array}
$$
Применяя к последнему равенству операцию усреднения и учитывая граничные условия, получим следующие четыре задачи
\begin{equation}\label{equ2}
\left\{\begin{array}{c}
\frac{\partial z^{(i)}_{-1}(\psi,\rho,\tau)}{\partial\tau}=\frac{\partial}{\partial\rho^2}z^{(i)}_{-1}(\psi,\rho,\tau),\\
\langle z^{(i)}_{-1}\rangle=0,\\
z^{(i)}_{-1}(\psi,\rho,\tau)|_{\rho=0}=-y^{(i)}_{-1},\\
z^{(i)}_{-1}(\psi,\rho,\tau)|_{\rho\to\infty}=0, i=1,2,
\end{array}\right.
\end{equation}
\begin{equation}\label{equ3}
\left\{\begin{array}{c}
\frac{\partial}{\partial\rho^2}v^{(i)}_{-1}(\psi,\rho)=0,\\
v^{(i)}_{-1}(\psi,\rho)|_{\rho\to\infty}=0, i=1,2,
\end{array}\right.
\end{equation}
\begin{equation}\label{equ4}
\left\{\begin{array}{c}
\frac{\partial}{\partial\rho}s_{-2}(\psi,\rho)=\frac{\partial}{\partial\tau}v^{(3)}_{-1}(\psi,\rho)=0,\\
s_{-2}(\psi,\rho)|_{\rho\to\infty}=0 
\end{array}\right.
\end{equation}
\begin{equation}\label{equ5}
\left\{\begin{array}{c}
\frac{\partial}{\partial\rho}n_{-2}(\psi,\rho,\tau)=\frac{\partial}{\partial\tau}z_{-1}^{(3)}(\psi,\rho,\tau)=0,\\
\langle n_{-2}\rangle = 0,\\
n_{-2}(\psi,\rho,\tau)|_{\rho\to\infty}=0. 
\end{array}\right.
\end{equation}
Из соотношений (\ref{equ2})-(\ref{equ5}) находим $z_{-1}^{(1,2)}\equiv 0$, $v_{-1}^{(1,2)}\equiv 0$, $s_{-2}\equiv 0$, $n_{-2}\equiv 0$.

Приравняем коэффициенты при степени $\omega$. Снова из уравнения $(\ref{eq:divC})$ получим
$$
\begin{array}{c}
\frac{\partial z_{0}^{(3)}(\psi,\rho,\tau)}{\partial\rho} = D_{-1,1}(\psi,\rho)z_{-1}^{(1)}(\psi,\rho,\tau)+D_{-1,2}(\psi,\rho)z_{-1}^{(2)}(\psi,\rho,\tau)+\\
+D_{-1,1}(\psi,\rho)z_{-1}^{(3)}(\psi,\rho,\tau),
\end{array}
$$
$$
\begin{array}{c}
\frac{\partial v_{0}^{(3)}(\psi,\rho)}{\partial\rho} = D_{-1,1}(\psi,\rho)v_{-1}^{(1)}(\psi,\rho)+D_{-1,2}(\psi,\rho)v_{-1}^{(2)}(\psi,\rho)+\\
+D_{-1,1}(\psi,\rho)v_{-1}^{(3)}(\psi,\rho).
\end{array}
$$
Откуда найдем $z_{0}^{(3)}(\psi,\rho,\tau)$ и $v_{0}^{(3)}(\psi,\rho)$. Для регулярных функций
\begin{equation}\label{equ6}
\left\{\begin{array}{c}
\frac{\partial y_{0}(x,\tau)}{\partial\tau}+\nabla m_{-2}(x,\tau)=\sum\limits_{1 \le |k| \le m}\sum\limits_{j=1}^{s}c^{j}_{-2}L_k(x)e_j(x) e^{ik\tau},\\
\Div y_{0}=0,\\
y_{0n}|_{\partial\Omega}=-z^{(3)}_{0}|_{\partial\Omega}|_{\rho=0}=0,\\
\langle y_0\rangle = \langle m_{-2}\rangle = 0.
\end{array}\right.
\end{equation}
Отсюда
\begin{equation}\label{equ7}
\left\{\begin{array}{c}
\Delta m_{-2}(x,\tau)=\Div \left(\frac{\partial y_0(x, \tau)}{\partial\tau}-\sum\limits_{1 \le |k| \le m} \sum\limits_{j=1}^{s}c^{j}_{-2}L_k(x)e_j(x)e^{ik\tau}\right)=0,\\
\frac{m_{-2}(x,\tau)}{\partial n}|_{\partial\Omega}=\left\{\frac{\partial y_0(x, \tau)}{\partial\tau}-\sum\limits_{1 \le |k| \le m} \sum\limits_{l=1}^{s}c^{l}_{-2}L_k(x)e_j(x)e^{ik\tau}\right\}_n,\\
\langle m_{-2}\rangle = 0,
\end{array}\right.
\end{equation}
где $n$ ~--- внутренняя нормаль к $\partial\Omega$. Теперь легко найдем $m_{-2}$ из задачи (\ref{equ7}) и $y_0$ из задачи (\ref{equ6}), зависящие от коэффициентов $c^j_{-2}$, $j=1,2,\dots,s$. При этом $y_{0}(x,\tau)=\sum\limits_{1 \le |k| \le m} \frac{L_k}{ik} \sum\limits_{j=1}^{s}c^{j}_{-2}e_j(x)e^{ik\tau}$. Из основного уравнения (\ref{eq:mainC}) для погранслойных функций получим
$$
\begin{array}{c}
\frac{\partial z_{0}(\psi,\rho,\tau)}{\partial\tau}=N_{-2}(\psi,\rho)[z_{0}(\psi,\rho,\tau)+v_{0}(\psi,\rho)]-\\
-P_{-1}(\psi,\rho)[s_{-1}(\psi,\rho)+n_{-1}(\psi,\rho,\tau)] + q_0(\psi,\rho)+ q_1(\psi,\rho,\tau),\\
\langle z_{0}\rangle = \langle n_{-1}\rangle = 0,
\end{array}
$$
где $q_0$, $q_1$ ~--- известные погранслойные вектор-функции, причем $q_1$ является $2\pi$-периодической с нулевым средним по $\tau$. Применяя операцию усреднения и учитывая граничные условия, получим задачи
\begin{equation}\label{eq:keyeq}
\left\{\begin{array}{c}
\frac{\partial z^{(i)}_{0}(\psi,\rho,\tau)}{\partial\tau}=\frac{\partial}{\partial\rho^2}z^{(i)}_{0}(\psi,\rho,\tau)+q_0(\psi,\rho),\\
\langle z^{(i)}_{0}\rangle=0,\\
z^{(i)}_{0}(\psi,\rho,\tau)|_{\rho=0}=-y^{(i)}_{0},\\
z^{(i)}_{0}(\psi,\rho,\tau)|_{\rho\to\infty}=0, i=1,2,
\end{array}\right.
\end{equation}
\begin{equation}
\left\{\begin{array}{c}
\frac{\partial}{\partial\rho^2}v^{(i)}_{0}(\psi,\rho)=q_1(\psi,\rho,\tau),\\
v^{(i)}_{0}(\psi,\rho)|_{\rho\to\infty}=0, i=1,2,
\end{array}\right.
\end{equation}
\begin{equation}
\left\{\begin{array}{c}
\frac{\partial}{\partial\rho}s_{-1}(\psi,\rho)=\frac{\partial}{\partial\tau}v^{(3)}_{0}(\psi,\rho)=0,\\
s_{-1}(\psi,\rho)|_{\rho\to\infty}=0 
\end{array}\right.
\end{equation}
\begin{equation}
\left\{\begin{array}{c}
\frac{\partial}{\partial\rho}n_{-1}(\psi,\rho,\tau)=\frac{\partial}{\partial\tau}z_{0}^{(3)}(\psi,\rho,\tau)=0,\\
\langle n_{-1}\rangle = 0,\\
n_{-1}(\psi,\rho,\tau)|_{\rho\to\infty}=0. 
\end{array}\right.
\end{equation}
Так найдем первые две компоненты вектора $z_0$, зависящие от неизвестных пока коэффициентов $c^j_{-2}$, $j=1,2,\dots,s$, а также $v_0^{(1,2)}$ и $s_{-1}$. 

Приравняем теперь коэффициенты при $\omega^{\frac{1}{2}}$. Для регулярных функций имеем
$$
\begin{array}{c}
\frac{\partial y_{1}(x,\tau)}{\partial\tau}+\nabla [p_{-1}(x)+m_{-1}(x,\tau)]=\\
=(\Delta+ B_0)u_{-1}(x)+\sum\limits_{1 \le |k| \le m}L_k(u_{-1}(x)+\sum\limits_{j=1}^{s}c^{j}_{-1}e_j(x))e^{ik\tau},
\end{array}
$$
применим операцию усреднения, получим задачу для $u_{-1}$
$$
\left\{\begin{array}{c}
(\Delta+ B_0)u_{-1}(x)=\nabla p_{-1}(x),\\
\Div u_{-1}(x)=0,\\
(u_{-1}(x),b_j(x))=0, j=1,2,\dots,s,\\
u_{-1}(x)|_{\partial\Omega}=v_{-1}|_{\rho=0}=0.
\end{array}\right.
$$
К задачам для последовательного нахождения $z^{3}_{1}$, $v^{3}_{1}$, $m_{-1}$, $y_{1}$, $z^{(1,2)}_{1}$, $v^{(1,2)}_{3}$, $s_0$ и $n_0$ вернемся после нахождения неизвестных коэффициентов $c^j_{-2}$, $j=1,2,\dots,s$.

Теперь приравняем коэффициенты при $\omega^0$
$$
\begin{array}{c}
\frac{\partial y_2(x,\tau)}{\partial\tau}+\nabla [p_0(x) + m_0(x,\tau)]=(\Delta+ B_0)(u_0(x)+y_0(x,\tau))+\\
+\sum\limits_{j=1}^{s}c^j_{-2}( C_0(x)+\sum\limits_{1\le |k|\le m}\frac{L_k(x)L_{-k}(x)}{ik})e_j(x)+\\
+\sum\limits_{1\le |k|\le m}(L_k(u_0(x) + \sum\limits_{j=1}^{s}c^{j}_{0}e_j(x)+y_0(x,\tau)) + d_k)e^{ik\tau}+ d_0(x),
\end{array}
$$
применим операцию усреднения
$$
\left\{\begin{array}{c}
-(\Delta+ B_0)u_0(x)+\nabla p_0(x)=\\
\sum\limits_{j=1}^{s}c^j_{-2}( C_0(x)+\sum\limits_{1\le |k|\le m}\frac{L_k(x)\Pi L_{-k}(x)}{ik})a_{j}(x)+ d_0(x),\\
\Div u_0(x) = 0,\\
(u_0(x),b_j(x))=0, j=1,2,\dots,s,\\
u_0(x)|_{\partial\Omega}=v_0(\psi,0).
\end{array}\right.
$$
Избавимся от неоднородности в граничном условии, подействуем проектором Вейля, и в качестве условия разрешимости для этой задачи получим систему линейныйх алгебраических уравнений
\begin{equation}\label{eq:SLAU1}
\begin{array}{c}
(\sum\limits_{j=1}^{s}c^j_{-2}(C_0(x)+\sum\limits_{1\le |k|\le m}\frac{L_k(x)L_{-k}(x)}{ik})a_{j}(x),b_i(x))=-(f_0(x),b_i(x)),\\
i=1,2,\dots,s,
\end{array}
\end{equation}
с известным столбцом свободных членов ($f_0(x)$ ~--- известная вектор-функция) и основной матрицей невырожденной матрицей $P$. Таким образом, система (\ref{eq:SLAU1}) однозначно разрешима. Возвращаясь к задаче (\ref{eq:keyeq}), однозначно находим $m_{-2}$, $y_{0}$, $z^{(1,2)}_0$. Теперь можно вернуться к нерешенным на предыдущем шаге задачам (см. выше). Продолжим этот процесс.

Допустим известны все коэффициенты вплоть до $y_{j-5}(x,\tau)$, $u_{j-5}(x)$, $z^{(1,2)}_{j-5}(\psi,\rho,\tau)$, $z^{3}_{j-4}(\psi,\rho,\tau)$, $v_{j-4}(\psi,\rho)$, $p_{j-5}(x)$, $m_{j-7}(x,\tau)$, $s_{j-4}(\psi,\rho)$, $n_{j-4}(\psi,\rho,\tau)$, $c^l_{j-5}$ включительно, $l=1,2,\dots,s$. Причем $y_{j-4}$, $m_{j-6}$, $z^{(1,2)}_{j-4}$ зависят от неизвестного пока коэффициента $c^l_{j-6}$, $l=1,2,\dots,s$. На $j$-м шаге, приравняв коэффициенты при $\omega^{-\frac{j-4}{2}}$, для регулярных функций получим
$$
\begin{array}{c}
\frac{\partial y_{j-3}(x,\tau)}{\partial\tau}+\nabla [p_{j-4}(x)+m_{j-4}(x,\tau)]=(\Delta + B_0)(u_{j-4}(x)+y_{j-4}(x,\tau)) +\\
+ C_0(x) (u_{j-6}(x)+y_{j-6}(x,\tau)) + \sum\limits_{1 \le |k| \le m}L_k(u_{j-4}(x)+y_{j-4}(x)+\\
+\sum\limits_{l=1}^{s}c^{l}_{j-4}a_l(x))e^{ik\tau},
\end{array}
$$
применим операцию усреднения, избавимся от неоднородности в граничных условиях и получим задачу для $u_{j-4}$
$$
\left\{\begin{array}{c}
-(\Delta+ B_0)u_{j-4}(x)+\nabla p_0(x)=\\
=\sum\limits_{l=1}^{s}c^l_{j-6}( C_0(x)+\sum\limits_{1\le |k|\le m}\frac{L_k(x)\Pi L_{-k}(x)}{ik})a_{l}(x)+ f(x),\\
\Div u_{j-4}(x) = 0,\\
(u_{j-4}(x),b_j(x))=0, j=1,2,\dots,s,\\
u_{j-4}(x)|_{\partial\Omega}=0.
\end{array}\right.
$$
Подействуем проектором Вейля, и в качестве условия разрешимости для этой задачи получим систему линейныйх алгебраических уравнений
\begin{equation}\label{eq:SLAU2}
\begin{array}{c}
(\sum\limits_{l=1}^{s}c^l_{j-6}(C_0(x)+\sum\limits_{1\le |k|\le m}\frac{L_k(x)L_{-k}(x)}{ik})a_{l}(x),b_i(x))=-(f(x),b_i(x)),\\
i=1,2,\dots,s,
\end{array}
\end{equation}
с известным столбцом свободных членов ($f(x)$ ~--- известная вектор-функция) и основной невырожденной матрицей вида $P$. Cистема (\ref{eq:SLAU2}) однозначно разрешима.Теперь однозначно находим $y_{j-4}$, $m_{j-6}$, $z^{(1,2)}_{j-4}$. Следующим шагом, как и показано выше, последовательно найдем $z^{(3)}_{j-3}$, $v^{(3)}_{j-3}$, $m_{j-5}$, $y_{j-3}$, $z^{(1,2)}_{j-3}$, $v^{(1,2)}_{j-3}$, $s_{j-4}$, $n_{j-4}$. Причем $m_{j-5}$, $y_{j-3}$, $z^{(1,2)}_{j-3}$ зависят от неизвестных пока коэффициентов $c^l_{j-5}$, $l=1,2,\dots,s$. Для погранслойных функций из основного уравнения (\ref{eq:mainC}) получим
\begin{equation}
\begin{array}{c}
\frac{\partial z_{j-2}(\psi,\rho,\tau)}{\partial\tau}+\sum\limits_{i=-1}^{j-4} P_{i}(\psi,\rho)[s_{j-4-i}(\psi,\rho)+n_{j-4-i}(\psi,\rho,\tau)]=\\
=\sum\limits_{i=-2}^{j-4}N_{i}(\psi,\rho)[z_{j-4-i}(\psi,\rho,\tau)+v_{j-4-1}(\psi,\rho)],\\
\langle z_{j-2}\rangle = \langle n_{j-3}\rangle = 0,
\end{array}
\end{equation}
т. е.
$$
\begin{array}{c}
\frac{\partial z_{j-2}(\psi,\rho,\tau)}{\partial\tau}=N_{-2}(\psi,\rho)[z_{j-2}(\psi,\rho,\tau)+v_{j-2}(\psi,\rho)]-\\
-P_{-1}(\psi,\rho)[s_{j-3}(\psi,\rho)+n_{j-3}(\psi,\rho,\tau)] + q_0(\psi,\rho,c^1_{j-5},c^2_{j-5},\dots,c^s_{j-5})+\\ +q_1(\psi,\rho,\tau,c^1_{j-5},c^2_{j-5},\dots,c^s_{j-5}),\\
\langle z_{j-2}\rangle = \langle n_{j-3}\rangle = 0,
\end{array}
$$
где $q_0$, $q_1$ ~--- известные вектор-функции, причем $q_1$ является $2\pi$-периодической с нулевым средним по $\tau$. Применяя операцию усреднения и учитывая граничные условия получим задачи
\begin{equation}
\left\{\begin{array}{c}
\frac{\partial z^{(i)}_{j-2}(\psi,\rho,\tau)}{\partial\tau}=\frac{\partial}{\partial\rho^2}z^{(i)}_{j-2}(\psi,\rho,\tau)+q_0(\psi,\rho,c^1_{j-5},c^2_{j-5},\dots,c^s_{j-5}),\\
\langle z^{(i)}_{j-2}\rangle=0,\\
z^{(i)}_{j-2}(\psi,\rho,\tau)|_{\rho=0}=-y^{(i)}_{j-2},\\
z^{(i)}_{j-2}(\psi,\rho,\tau)|_{\rho\to\infty}=0, i=1,2,
\end{array}\right.
\end{equation}
\begin{equation}
\left\{\begin{array}{c}
\frac{\partial}{\partial\rho^2}v^{(i)}_{j-2}(\psi,\rho)=q_1(\psi,\rho,\tau,c^1_{j-5},c^2_{j-5},\dots,c^s_{j-5}),\\
v^{(i)}_{j-2}(\psi,\rho)|_{\rho\to\infty}=0, i=1,2,
\end{array}\right.
\end{equation}
\begin{equation}
\left\{\begin{array}{c}
\frac{\partial}{\partial\rho}s_{j-3}(\psi,\rho)=\frac{\partial}{\partial\tau}v^{(3)}_{j-2}(\psi,\rho)=0,\\
s_{j-3}(\psi,\rho)|_{\rho\to\infty}=0 
\end{array}\right.
\end{equation}
\begin{equation}
\left\{\begin{array}{c}
\frac{\partial}{\partial\rho}n_{j-3}(\psi,\rho,\tau)=\frac{\partial}{\partial\tau}z_{j-2}^{(3)}(\psi,\rho,\tau)=0,\\
\langle n_{j-3}\rangle = 0,\\
n_{j-3}(\psi,\rho,\tau)|_{\rho\to\infty}=0. 
\end{array}\right.
\end{equation}
К этим задачам вернемся после нахождения неизвестных коэффициентов $c^l_{j-5}$, $l=1,2,\dots,s$.

\subsection{Обоснование оценок асимптотики}

В заключение докажем оценки (\ref{eq:ev}). В силу построения асимптотики частичная сумма $u_\omega^j$ удовлетворяет уравнению в операторной форме
\begin{equation}\label{eq:oper1}
\frac{\partial u}{\partial t} = \left(A + \frac{1}{\omega}K\right)u + \sum\limits_{1\le |k|\le m}\left(M_k(x)u+a_k(x)\right)e^{ik\omega t}+a_0(x) + r_\omega^j(x,t),
\end{equation}
где при любом $l\ge 0$
\begin{equation}\label{eq:ost1}
\|r_\omega^j\|_{C^{l,\frac{l}{2}}(\bar{Q})}\le c_j\omega^{-\frac{j+1-l}{2}}.
\end{equation}
Из (\ref{eq:oper}), (\ref{eq:oper1}) следует, что разность $w_\omega^j=u_\omega-u_\omega^j$ является решением уравнения
\begin{equation}\label{eq:oper2}
\frac{\partial w}{\partial t} - \left(A + \frac{1}{\omega}K\right)w - \sum\limits_{1\le |k|\le m}M_k(x)w e^{ik\omega t} = - r_\omega^j(x,t).
\end{equation}
Из оценок (\ref{eq:ost}), (\ref{eq:ost1}) следует неравенство
\begin{equation}\label{eq:ost2}
\|w_\omega^j\|_{C(R,L_2(\Omega))}\le d_j\omega^{-\frac{j+1}{2}}.
\end{equation}
Очевидно, вектор-функция $v_\omega^j=\chi(t)w_\omega^j$, $t\in[0,1]$ (определение $\chi(t)$ см. выше) является решением однородной начально-краевой задачи для уравнения
\begin{equation}\label{eq:koshi1}
\begin{array}{c}
\frac{dv}{dt}-Av=\chi\left[\left(\frac{1}{\omega}K + \Pi\sum\limits_{1\le |k|\le m}L_k(x)e^{ik\omega t}\right)v - r_\omega^j(x,t) \right]-\frac{d\chi}{dt}v
\end{array}
\end{equation}
Многократно применяя к (\ref{eq:koshi1}) теорему 2 \cite{Solon}, с учетом (\ref{eq:ost2}), выводим оценки (\ref{eq:ev}). Теорема 3 доказана.

\end{document}